\documentclass[a4paper,10pt]{article}

\usepackage{amsthm}

\usepackage{amsmath}

\usepackage{amsfonts}

\usepackage{amssymb}

\usepackage{amsmath,amsthm,amscd,amssymb}

\usepackage{latexsym}

\usepackage{graphicx}

\usepackage[all]{xy}

\theoremstyle{plain}

\newtheorem{thm}{\textbf{Theorem}}

\newtheorem{lem}{\textbf{Lemma}}

\newtheorem{df}{\textbf{Definition}}

\newtheorem{cor}{\textbf{Corollary}}

\newtheorem{prop}{\textbf{Proposition}}

\newcommand{\A}{\Bbb{A}}

\newcommand{\R}{\Bbb{R}}

\newcommand{\C}{\Bbb{C}}

\newcommand{\HH}{\mathcal{H}}

\newcommand{\Q}{\Bbb{Q}}

\newcommand{\G}{\Bbb{G}}

\newcommand{\Z}{\Bbb{Z}}

\newcommand{\T}{\Bbb{T}}

\newcommand{\p}{\frak{p}}

\newcommand{\q}{\frak{q}}

\newcommand{\g}{\frak{g}}

\newcommand{\vv}{\tilde{v}}

\newcommand{\VV}{\Bbb{V}}

\newcommand{\m}{\frak{m}}

\newcommand{\K}{\mathcal{K}}

\newcommand{\CC}{\mathcal{C}}

\newcommand{\OO}{\mathcal{O}}

\newcommand{\LL}{\mathcal{L}}

\newcommand{\XX}{\mathcal{X}}

\newcommand{\RR}{\mathcal{R}}

\newcommand{\Frob}{\text{Frob}}

\newcommand{\Art}{\text{Art}}

\newcommand{\Rep}{\text{Rep}}

\newcommand{\Res}{\text{Res}}

\newcommand{\Hom}{\text{Hom}}

\newcommand{\Ext}{\text{Ext}}

\newcommand{\End}{\text{End}}

\newcommand{\tr}{\text{tr}}

\newcommand{\Gal}{\text{Gal}}

\newcommand{\GL}{\text{GL}}

\newcommand{\Spec}{\text{Spec}}

\newcommand{\SL}{\text{SL}}

\newcommand{\y}{\hspace{6pt}}

\title{{\bf{Weak local-global compatibility in the $p$-adic Langlands program for $U(2)$}}}

\author{Przemyslaw Chojecki $\y$ Claus Sorensen}

\begin{document}

\date{}

\maketitle

\begin{abstract}
We initiate a careful study of the completed cohomology $\hat{H}^0$ of the tower of arithmetic manifolds of $G$, a definite unitary group in two variables associated with a CM-extension $\K/F$. When the prime $p$ splits, we prove that (under technical asumptions) the $p$-adic local Langlands correspondence for $\GL_2(\Q_p)$ occurs in $\hat{H}^0$, strongly inspired by Emerton's work on local-global compatibility for $\GL(2)_{/\Q}$. As an application, we obtain a result towards the Fontaine-Mazur conjecture over $\K$
(roughly saying that a two-dimensional geometric Galois representation arising from a $p$-adic modular form must in fact be modular). In fact we prove something slightly stronger. If $x$ is a point on the eigenvariety such that $\rho_x$ is geometric (and satisfying additional hypotheses which we suppress), then $x$ must be a classical point.
Thus, not only is $\rho_x$ modular, but the weight of $x$ defines an accessible refinement. This follows from a recent result of Colmez (which describes the locally analytic vectors in $p$-adic unitary principal series), knowing that $\rho_x$ admits a triangulation compatible with the weight (as shown recently by Hellmann).

\footnote{{\it{2000 AMS Mathematics Classification}}: 11F33.}

\end{abstract}


\section{Introduction}

In the last decade, the $p$-adic Langlands program has exploded with activity. In a nutshell, it predicts a close relationship between continuous representations
$\Gamma_{F} = \Gal (\bar{F} / F) \rightarrow \GL_n(E)$, where $F$ and $E$ are both finite extensions of $\Q_p$,  and unitary representations of $\GL_n(F)$ on Banach $E$-spaces.
For $\GL_2(\Q_p)$, this correspondence was pioneered by Breuil, and in the last few years one has achieved a more complete understanding for $\GL_2(\Q_p)$, due to the work of many people (notably Berger, Colmez, Emerton, Kisin, Paskunas, and others). These results have already had astounding applications to various notoriously difficult problems in number theory.

\medskip

\noindent This article takes its point of departure in Emerton's progress on the celebrated Fontaine-Mazur conjecture (which roughly gives a local characterization of Galois representations arising from geometry). In \cite{Em1}, and its predecessor \cite{Em0}, he explains how the $p$-adic local Langlands correspondence for $\GL_2(\Q_p)$ appears in the completed cohomology $\hat{H}^1$ of the tower of modular curves. From the "de Rham" condition, one gets the existence of locally algebraic vectors, and Emerton deduces that a "promodular" representation $\Gamma_{\Q} = \Gal (\bar{\Q} / \Q) \rightarrow \GL_2(E)$, which is de Rham (with distinct Hodge-Tate weights), must in fact be modular; again, under weak technical assumptions, which we will not record here. 

\medskip

\noindent In this paper we look at the tower of arithmetic manifolds of a definite unitary group in two variables, $G=\Res_{F/\Q}(U)$; an inner form of $\GL(2)_{/\K}$ over a CM extension $\K/F$. The arithmetic manifolds are in fact just (arithmetically rich) finite sets, occasionally called "Hida varieties". Thus, instead of $\hat{H}^1$, we are looking at 
the Banach space $\hat{H}^0$, which can be realized as the space of continuous functions on a profinite set (endowed with the sup-norm). When the prime $p$ splits, we relate the $p$-adic local Langlands correspondence for $\GL_2(\Q_p)$ to $\hat{H}^0$ (see Theorem 1 below for a precise statement). As a result thereof, we obtain a corollary towards the Fontaine-Mazur conjecture for representations $\Gamma_{\K}\rightarrow \GL_2(E)$ associated with $p$-adic modular forms on $G$, in the vein of Emerton (see Theorem 2 below). In fact we prove a little more than modularity; we prove classicality (that is, we keep track of refinements and triangulations). 

\medskip

\noindent To orient the reader, we briefly point out the major differences with \cite{Em1}:

\begin{itemize}
\item We work with $\hat{H}^0$, which we find to be more hands-on than $\hat{H}^1$. The arithmetic manifolds of $G$ are not Shimura varieties; so $\hat{H}^0$ carries no Galois-action (as opposed to $\hat{H}^1$ of modular curves). This simplifies some of the arguments. 

\item We allow $F \neq \Q$, but must assume $p$ splits in $\K$. Thus we really embed a tensor product (over places $v|p$ of $F$) of 
$p$-adic local Langlands correspondents in $\hat{H}^0$. Consequently, we make progress on Fontaine-Mazur for two-dimensional representations of $\Gamma_{\K}$, for
CM-fields $\K$, as opposed to $\Gamma_{\Q}$.

\item We work with a fixed tame level $K^p$ throughout, to surmount the difficulty with the non-uniqueness of hyperspecials in $p$-adic $U(2)$; 
and find it convenient to formulate our results in terms of the eigenvariety $X=X_{K^p}$. 

\item We prove classicality, not "just" modularity: That is, if the representation $\rho_x$ at $x \in X$ is de Rham (with distinct Hodge-Tate weights etc.) then $x$ is a classical point. Here the crux of the matter lies in relating weights and refinements (which makes critical use of recent work of Colmez, and of Hellmann's results on triangulinity). 

\item On the flip-side, for now, we must make the rather bold assumption that the mod $p$ reduction $\bar{\rho}_x$ is irreducible at all places of $\K$ above $p$. Emerton gets by with much weaker hypotheses at $p$, using \cite{BE}.
\end{itemize}

\noindent In order to state our main results, we must briefly set up the notation in use throughout the paper. Once and for all, we fix a prime number $p$. To be safe, we will always assume $p>3$. We let $F/\Q$ be a totally real field, and $\K/F$ a CM extension, in which $p$ {\it{splits}} completely. Places of $F$ are usually denoted by $v$, and those of $\K$ by $w$. 
For each place $v|p$ of $F$, we choose a place $\vv|v$ of $\K$ above it (note that $\K_{\vv}=\Q_p$, canonically). Given an algebraic isomorphism $\iota: \C \longrightarrow \bar{\Q}_p$, the choice of a collection $\{\vv\}$ amounts to choosing a CM-type, which is ordinary for $\iota$, in the sense of Katz. We will use it to identify a certain $p$-adic group 
$G(\Q_p)$ with a product of copies of $\GL_2$, in a compatible way.

\medskip

\noindent Let $D$ be a quaternion algebra over $\K$, endowed with an $F$-linear anti-involution $\star$ of the second kind ($\star|_{\K}=c$). This pair defines a unitary group $U=U(D,\star)_{/F}$, an inner form of $\GL(2)$ over $\K$. Indeed, $U \times_F \K \simeq D^{\times}$. We find it convenient to work over $\Q$, and introduce 
$G=\text{Res}_{F/\Q}(U)$. We will always assume $G(\R)$ is compact, and that $D$ splits above $p$. Using our choices $\{\vv\}$, we identify
$$
\text{$G(\R)\overset{\sim}{\longrightarrow} U(2)^{\Hom(F,\R)}$, $\y$ $G(\Q_p) \overset{\sim}{\longrightarrow} \prod_{v|p} \GL_2(\K_{\vv})$.}
$$
(Of course, $\K_{\vv}=\Q_p$, but we wish to incorporate $\vv$ in our notation to emphasize how our identification depends on this choice. Hence we stick to the somewhat cumbersome notation $\K_{\vv}$.)

\medskip

\noindent Our main occupation will be Galois representations $\rho$ associated with $p$-adic modular forms on $G$. We will show a "Fontaine-Mazur" like result of the following form: If $\rho$ is geometric, with distinct weights (and satisfies additional technical hypotheses), then in fact $\rho$ arises from a {\it{classical}} modular form on $G$.

\medskip

\noindent As is well-known, $p$-adic modular forms (of finite slope), and their Galois representations, are nicely parametrized by rigid analytic spaces called eigenvarieties. For $U(n)$, they have been constructed in great detail in \cite{BC} by Bellaiche-Chenevier (for $F=\Q$), and in \cite{Che} by Chenevier (for any $F$), elaborating on his thesis,
and using Buzzard's "eigenvariety machine" from \cite{Buz}. A more general construction was given by Emerton in \cite{Em} (which certainly covers the case where $G(\R)$ is compact; and it does not require Iwahori-level at $p$). Yet another construction, in the style of Chenevier, was given by Loeffler in \cite{Loe}, only assuming $G(\R)$ is compact (and curiously dealing with any parabolic, as opposed to just the Borel). By the uniqueness, shown in \cite{BC} for example, all these constructions are compatible, and define the same eigenvariety. In the special case where $F=\Q$, the eigenvariety for our two-variable unitary group $G$, is formally reminiscent of the mother of all eigenvarieties; the celebrated "eigencurve" of Coleman and Mazur \cite{CM}, which in turn has its origin in Hida theory (the slope zero case).

\medskip

\noindent We will not use much about eigenvarieties, besides their definition and basic structural properties, but we find them to provide an elegant framework for our results.
It is pleasing for us to work with a fixed tame level $K^p=\prod_{v\nmid p}K_v$ throughout the paper. Given $K^p$, let $\Sigma_0$ denote the set of places $v \nmid p$ for which $K_v$ is not a hyperspecial subgroup, and then introduce $\Sigma=\Sigma_0\sqcup \Sigma_p$. The eigenvariety $X=X_{K^p}$ parametrizes eigensystems of $\HH(K^p)^{sph}$ (the Hecke algebra of $G(\A_f^{\Sigma})$ relative to $K^{\Sigma}$), which are associated with $p$-adic modular forms. To be more precise, we first introduce weight space
$$
\text{$\hat{T}=\Hom_{loc.an.}(T(\Q_p),\G_m^{rig})$, $\y$ $T(\Q_p)\simeq \prod_{v|p}T_{\GL(2)}(\K_{\vv})$.}
$$
This is a rigid analytic space over $\Q_p$, endowed with a universal locally analytic character $\hat{T}\rightarrow \OO(\hat{T})^{\times}$. See paragraph 8 of \cite{Buz} for more details. Similarly, $\mathcal{W}$ represents locally analytic characters of $T(\Z_p)$, a disjoint union of open balls. Moreover, $\hat{T}$ is non-canonically isomorphic to the direct product $\mathcal{W} \times (\G_m^{rig})^d$.

\medskip

\noindent Now, the eigenvariety $X=X_{K^p}$ is a reduced rigid analytic variety over $\Q_p$ (since $\K_{\vv}=\Q_p$), which comes equipped with a finite morphism 
$\chi: X \rightarrow \hat{T}$, and additional structure (see Theorem 1.6 on p. 5 in \cite{Che}, for example):

\begin{itemize}
\item $\lambda: \HH(K^p)^{sph} \rightarrow \OO(X)$, a homomorphism of $\Q_p$-algebras,
\item $t: \Gamma_{\K}\rightarrow \OO(X)$, a pseudo-character,
\item $X_{cl}\subset X(\bar{\Q}_p)$, a Zariski-dense subset.
\end{itemize}

\noindent Here $t$ is "associated" with $\lambda$, in a natural sense (see 2.1 below). Furthermore,
$$
\text{$X(\bar{\Q}_p)\longrightarrow (\hat{T} \times \Spec  \HH(K^p)^{sph})(\bar{\Q}_p)$, $\y$ $x \mapsto (\chi_x,\lambda_x)$,}
$$
restricts to a bijection between $X_{cl}$ and the set of "classical" points. That is, those $(\chi,\lambda)$, for which $\chi=\psi\theta$ is locally algebraic, and there exists an automorphic representation $\pi$ of $G(\A)$, of weight $\psi$ and tame level $K^p$, such that $\pi_p$ embeds into the (unnormalized) principal series representation $i_B(\theta)$. 

\medskip

\noindent Finally, for each $x \in X(\bar{\Q}_p)$, we let $\rho_x: \Gamma_{\K}\rightarrow \GL_2(\bar{\Q}_p)$ denote the unique semisimple continuous representation with 
$\tr \rho_x=t_x$. When $x \in X_{cl}$, it is known that $\rho_x$ is "geometric" in the sense of Fontaine-Mazur (unramified almost everywhere, and de Rham above $p$). Our goal 
here is a partial converse. We first state the weak local-global compatibility, referred to in the title, which roughly says the $p$-adic local Langlands correspondence occurs in the completed cohomology of the tower of arithmetic manifolds of $G$ (as announced on p. 8 in \cite{Sor}). We are hoping to write a sequel proving strong local-global compatibility, which gives a more precise description of the completed cohomology, which involves "local Langlands in families" away from $p$ (following Emerton-Helm \cite{EH}). However, for the applications to Fontaine-Mazur, the weak version suffices. 

\begin{thm}
Let $E/\Q_p$ be a finite extension, and let $x \in X(E)$ be a point on the eigenvariety, of tame level $K^p$, which satisfies the following assumptions: 
\begin{itemize}
\item[(1)] $\rho_{x}$ is defined over $E$.
\item[(2)] $\bar{\rho}_{x,w}$ is absolutely \underline{irreducible}, for all $w|p$.
\end{itemize}
Then there is a nonzero continuous $G(\Q_p)$-equivariant map ($\p_x=\ker \lambda_x$),
$$
\hat{\otimes}_{v|p} B(\rho_{x,\vv})\longrightarrow \hat{H}^0(K^p)_E[\p_x],
$$
where $B(\cdot)$ is the $p$-adic local Langlands correspondence for $\GL_2(\Q_p)$.
\end{thm}

\medskip

\noindent The flaring assumption here is (2); the mod $p$ irreducibility at $p$. This simplifies many of the arguments a great deal, and we currently do not know how to relax this condition. It would require an analogue of \cite{BE}, which is based on a geometrically intricate construction of "overconvergent" companion forms. It is not clear to us how (and if) this can be done in our unitary group context, but it is something we hope to work on. Guided by \cite{Em1}, one would expect that (2) can be replaced with the following two (much weaker) assumptions:

\begin{itemize}
\item[(3)] $\bar{\rho}_x$ is absolutely irreducible (as a representation of the full $\Gamma_{\K}$).
\item[(4)] $\bar{\rho}_{x,w} \nsim \chi \otimes \begin{pmatrix}1 & * \\ & \bar{\epsilon}\end{pmatrix}$, for all places $w|p$ of $\K$.
\end{itemize}

\noindent We wish to stress that the proof of Theorem 1 hinges upon local-global compatibility at $p=\ell$, now known (for $U(n)$ even) due to the hard work of Barnet-Lamb, Gee, Geraghty, and Taylor (see \cite{B1} and its follow-up \cite{B2}), and that of Caraiani \cite{Car}.

\medskip

\noindent Our main motivation for proving this weak local-global compatibility, is its application to the Fontaine-Mazur conjecture, along the following lines:

\begin{thm}
Let $x \in X(E)$ be a point as in Theorem 1, and assume moreover that $\rho_{x,\vv}$ is potentially semistable, with distinct Hodge-Tate weights, for all $v|p$.
Then $x$ is a classical point, $x \in X_{cl}$. 
\end{thm}

\medskip

\noindent This is more than a modularity result. To conclude $x\in X_{cl}$, one needs to keep careful track of refinements. This involves the existence of specific triangulations
of $\rho_{x,\vv}$, proved in many instances by Hellmann, and a description of the locally analytic vectors in unitary principal series; a recent result due to Colmez \cite{Co}.

\medskip

\noindent {\it{Acknowledgments}}. The tremendous debt this paper owes to Emerton's work \cite{Em1} will be clear to the experts. We wish to acknowledge its impact on both of us, as well as that of Breuil's beautiful Bourbaki survey \cite{Bre}, which helped us getting into this circle of ideas. Finally, we wish to thank the Fields Institute in Toronto for hosting the $p$-adic Langlands workshop in April, 2012, where this collaboration begun.

\section{Automorphic Galois representations}

Galois representations associated with regular polarized cusp forms on $\GL(n)_{/K}$ are now almost completely understood, thanks to the "book project" and its spin-offs. This is the culmination of collective efforts of many people, initiated by Clozel, Harris, Kottwitz, Taylor, and others. Here we wish to briefly give the lay of the land for unitary groups in two variables. 

\medskip

\noindent Thus let $\pi$ be an automorphic representation of $G(\A)$, whose infinity component $\pi_{\infty}=\otimes_{v|\infty}\pi_v$ is an irreducible algebraic representation of $G(\C)$, restricted to the compact Lie group $G(\R)$. We assign weights to $\pi_{\infty}$ as follows: Each factor $\pi_v$ is a representation of $U(F_v)$, which we may "complexify" to a representation of $\GL_2(\K_w)$, where $w|v$ is unique. Upon choosing an embedding $\tau: \K \hookrightarrow \C$, among the pair corresponding to $w$, we identify this with a representation of $\GL_2(\C)$, which is irreducible algebraic of highest weight $\psi_{\tau}$, relative to the lower-triangular Borel. In other words, there are integers 
$\kappa_{1,\tau}< \kappa_{2,\tau}$ such that
$$
\text{$\psi_{\tau}(t)=t_1^{\kappa_{1,\tau}}t_2^{\kappa_{2,\tau}-1}$, $\y$ $t = \left(\begin{smallmatrix} t_1 &  \\  & t_2 \end{smallmatrix} \right)  \in T_{\GL(2)}(\C)$.}
$$
(Note that $\kappa_{1,\bar{\tau}}=1-\kappa_{2,\tau}$ and $\kappa_{2,\bar{\tau}}=1-\kappa_{1,\tau}$.) Another way to think of this is in terms of the local base change $BC_{w|v}(\pi_v)$, which has infinitesimal character given by $\psi_{\tau}$, upon identifying it with a representation of $\GL_2(\C)$, via $\tau$.

\begin{thm}
Choose an isomorphism $\iota: \C \overset{\sim}{\longrightarrow} \bar{\Q}_p$. With $\pi$ as above, we may associate a unique continuous semisimple Galois representation,
$$
\rho=\rho_{\pi,\iota}: \Gamma_{\K}=\Gal(\bar{\Q}/\K) \rightarrow \GL_2(\bar{\Q}_p),
$$
satisfying the following list of desiderata:
\begin{itemize}
\item For every finite place $w|v$, for which the local base change $BC_{w|v}(\pi_v)$ is defined (even those above $p$), one has
$$
WD(\rho|_{\Gamma_{\K_w}})^{F-ss}\simeq \iota rec(BC_{w|v}(\pi_v)\otimes |\det|^{-1/2}).
$$
(Here the classical local Langlands correspondence $rec(\cdot)$ is normalized as in \cite{HT}. At $w|p$, the Weil-Deligne representation is defined by Fontaine.)
\item Indeed, for every place $w|p$, the restriction $\rho_w=\rho|_{\Gamma_{\K_w}}$ is potentially semistable, with labelled Hodge-Tate weights determined by $\pi_{\infty}$,
$$
HT_{\tau}(\rho_w)=\{\kappa_{1,\tau}< \kappa_{2,\tau}\},
$$
for $\tau: \K_w \hookrightarrow \bar{\Q}_p$, tacitly identified with an embedding $\K \hookrightarrow \C$ via $\iota$. 
\item $\rho^{\vee}\simeq \rho^c \otimes \epsilon$ (where $\epsilon=\epsilon_{cyc}$ is the cyclotomic character of $\Gamma_{\K}$).
\item $\det \rho \circ \Art_{\K}=BC_{\K/F}(\chi_{\pi})\cdot \epsilon^{-1}$.
\end{itemize}
(The local base change $BC_{w|v}(\pi_v)$ is defined when $v$ splits in $\K$, or $\pi_v$ is unramified; for some hyperspecial maximal subgroup.) 
\end{thm}

\noindent {\it{Proof}}. Let us quickly sketch the argument, and list the key references. Using Rogawski's book \cite{Rog}, we first base change $\pi$ to $\GL(2)_{/\K}$, resulting in an "isobaric" automorphic representation $\Pi=BC_{\K/F}(\pi)$, whose infinity component $\Pi_{\infty}$ has the same infinitesimal character as an algebraic representation. Moreover,
$\Pi^{\vee}\simeq \Pi^c$. If $\Pi$ is cuspidal, we take $\rho=\rho_{\Pi,\iota}$, where the latter is given by Theorem 1.2 in \cite{CY}, for instance. It satisfies the desired properties by 
Theorem 1.1 in \cite{BLGGT}, and Caraiani's sequel \cite{Car}, which removes the Shin-regularity assumption. At last, if $\Pi$ is non-cuspidal, it is an isobaric sum of algebraic Hecke characters $\chi_1 \boxplus \chi_2$, with which we associate Galois characters $\rho_{\chi_i,\iota}$ via class field theory. Then we let
$\rho=\rho_{\chi_1,\iota}\oplus \rho_{\chi_2,\iota}$, for which it is straightforward to verify the properties. $\square$

\medskip

\noindent {\it{Remark}}. The analogous result holds for unitary groups in $n$ variables, except that a certain "regularity" condition creeps in (ruling out $\pi_{\infty}=1$ for example), which ensures that the base change is an isobaric sum of cusp forms (as opposed to just discrete automorphic representations). The lucky circumstance for $n=2$, which we rely on, is that "discrete" is the same as "cuspidal" for $\GL(1)$.

\medskip

\noindent For later use, we spell out what happens at the places $v$ where $\pi$ is unramified. First off, $U_{/F_v}$ must be unramified (quasi-split and split over an unramified extension), in which case it has two conjugacy-classes of hyperspecial subgroups. Thus we should really pick one, say $K_v$, and specify that $\pi_v^{K_v}\neq 0$. Via the Satake parametrization, local base change is given by pulling back eigensystems along the natural algebra-homomorphism (defined over $\Z[p^{\frac{1}{2}}]$),
$$
b_{w|v}: \HH(\GL_2(\K_w),\GL_2(\OO_w)) \twoheadrightarrow \HH(U(F_v),K_v).
$$
For the surjectivity of this map, and a thorough discussion of the salient facts, see \cite{Min} (especially Corollary 4.2). For each $w$, we let $T_w$ be the Hecke operator 
for $\GL_2(\K_w)$, which acts on $\Pi_w^{\GL_2(\OO_w)}$ by the sum of the integral Satake parameters of $\Pi_w$. Similarly, $S_w$ acts by their product. Then, the target Hecke algebra $\HH(K_v)$ is generated by $t_w=b_{w|v}(T_w)$ and $s_w=b_{w|v}(S_w)$. Furthermore,
$$
\text{$\tr \rho_{\pi,\iota}(\Frob_w)=\iota \lambda_{\pi_v}(t_w)$, $\y$ $\det \rho_{\pi,\iota}(\Frob_w)=\iota \lambda_{\pi_v}(s_w)$,}
$$
where $\lambda_{\pi_v}: \HH(K_v)\rightarrow \C$ is the eigensystem of $\pi_v^{K_v}$, when $\pi_v$ is $K_v$-unramified.


\section{Theorem 1 $\Rightarrow$ Theorem 2}

\subsubsection{Classical weight}

First off, let us note that $x$ at least has a classical weight $\chi_x=\psi_x\theta_x$, where $\theta_x$ is smooth, and $\psi_x$ is a dominant algebraic character
of $T(\Q_p)$. Indeed, for any $x$, the Hodge-Tate-Sen weights $\{\kappa_1,\kappa_2\}$ of $\rho_{x,\vv}$ are encoded in $\chi_x=\otimes_{v|p}\chi_{x,\vv}$ as follows. On a compact open subgroup of $T_{\GL(2)}(\K_{\vv})$, the character $\chi_{x,\vv}$ takes the form
$$
\chi_{x,\vv}(t)=t_1^{\kappa_1}t_2^{\kappa_2-1}.
$$
Since $\rho_{x,\vv}$ is Hodge-Tate, $\kappa_1<\kappa_2$ are integers, and $\psi_{x,\vv}(t)=t_1^{\kappa_1}t_2^{\kappa_2-1}$ is an algebraic character on all of 
$T_{\GL(2)}(\K_{\vv})$, which is dominant (relative to the lower-triangular Borel), and agrees with $\chi_{x,\vv}$ in a neighborhood of the identity.

\subsubsection{Modularity}

\noindent By assumption, $\rho_{x,\vv}$ is absolutely irreducible. By well-known properties of the $p$-adic local Langlands correspondence, we infer that $B(\rho_{x,\vv})$ is topologically irreducible. Thus the locally algebraic vectors $B(\rho_{x,\vv})^{alg}$ are dense, being nonzero by assumption (1) of Theorem 2. Consequently, $\otimes_{v|p}B(\rho_{x,\vv})^{alg}$
is dense in $\hat{\otimes}_{v|p} B(\rho_{x,\vv})$, and by continuity of the map in Theorem 1, it restricts to a nonzero map, equivariant under the $G(\Q_p)$-action,
$$
\otimes_{v|p}B(\rho_{x,\vv})^{alg} \longrightarrow \hat{H}^0(K^p)_E[\p_x]^{alg}.
$$
In particular, the target is non-trivial, which is to say there is an automorphic representation $\pi$ of $G(\A)$, with coefficients in $E$, such that $\pi_{\infty}$ has highest weight $\psi_x$, and the spherical Hecke algebra $\HH(K^p)^{sph}$ acts on $\pi_f^{K^p}\neq 0$ via the eigensystem $\lambda_x: \HH(K^p)^{sph} \rightarrow E$, with kernel $\p_x$.
This shows that $\rho_x$ is at least modular; it is associated with the automorphic representation $\pi$. To conclude that $x \in X_{cl}$, we are left with verifying that there is an embedding $\pi_p \hookrightarrow i_B(\theta_x)$. Here $B$ is the product of the upper-triangular 
$B_{\GL(2)}(\K_{\vv})$, and the induction is not unitarily normalized (no modulus factor $\delta_B^{\frac{1}{2}}$ involved). We say $\theta_x$ is an accessible refinement of 
$\pi_p=\otimes_{v|p}\pi_{v}$. 

\subsubsection{Refinements}

We will use the following result, which is work in progress, of Hellmann. 

\begin{thm}
When $\chi_{x,\vv}$ is regular, $\rho_{x,\vv}$ is trianguline. In fact, its \'etale $(\phi,\Gamma)$-module $D_{rig}(\rho_{x,\vv})$ admits a triangulation,
$$
0 \longrightarrow \RR(\delta_1) \longrightarrow D_{rig}(\rho_{x,\vv}) \longrightarrow \RR(\delta_2) \longrightarrow 0,
$$
where $\RR$ is the Robba ring over $E$, and the $\delta_i$ are $E$-valued (necessarily) continuous characters of $\Q_p^{\times}=\K_{\vv}^{\times}$, related to the weight $\chi_{x,\vv}$ via the formula
$$
\chi_{x,\vv}=\delta_1 \otimes \delta_2 \epsilon^{-1},
$$
where $\epsilon$ is the "cyclotomic" character of $\Q_p^{\times}$, which kills $p$ (sending $x \mapsto x|x|$).
\end{thm}

\noindent {\it{Proof}}. When $F=\Q$, this is Corollary 4.3 in \cite{He}, which in fact gives an analogue for $\GL(n)$. (This uses the residual irreducibility of $\rho_x$.) In the current version of \cite{He}, the normalization of weights is a bit off (the $\epsilon$-shift between the highest weights and the Hodge-Tate-Sen weights is ignored), but this can easily be fixed. Hellmann has informed us he is writing up a sequel, for any $F$, in which this will be corrected. $\square$ 

\medskip

\noindent {\it{Remark}}. We will not be precise about what we mean by $\chi_{x,\vv}$ being regular (see \cite{He} page 7 for the definition). It is automatically satisfied if $\rho_{x,\vv}$ is potentially semistable with distinct Hodge-Tate weights $\kappa_1<\kappa_2$, which we assume.

\medskip

\noindent Consequently, $B(\rho_{x,\vv})$ is a $p$-adic unitary principal series for $\GL_2(\Q_p)$, and Colmez has recently described its locally analytic vectors 
$B(\rho_{x,\vv})^{an}$, thereby proving conjectures of Breuil and Emerton. A different proof was given by Liu, Xie, and Zhang. We recall Colmez's result below.

\begin{thm}
Suppose $\delta_1\delta_2^{-1}\neq x^k|x|$, for any $k \in \Z_+$, and $\rho_{x,\vv}$ is irreducible. Then $B(\rho_{x,\vv})^{an}$ sits in an exact sequence of locally analytic representations,
$$
0 \longrightarrow i_B(\delta_2\otimes \delta_1\epsilon^{-1})^{an}\longrightarrow B(\rho_{x,\vv})^{an} \longrightarrow i_B(\delta_1\otimes \delta_2\epsilon^{-1})^{an}\longrightarrow 0.
$$
(Here the induction is unnormalized.)
\end{thm}

\noindent {\it{Proof}}. This is Theorem 0.7, part (i), on page 7 in \cite{Co}. See also the main Theorem 1.2 in \cite{LXZ}, which gives an alternative proof. $\square$

\medskip

\noindent Taking locally algebraic vectors, which is left exact, yields 
$$
0 \longrightarrow i_B(\delta_2\otimes \delta_1\epsilon^{-1})^{alg}\longrightarrow B(\rho_{x,\vv})^{alg} \longrightarrow i_B(\delta_1\otimes \delta_2\epsilon^{-1})^{alg}.
$$
Here the rightmost term is $i_B(\chi_{x,\vv})^{alg}$, by Theorem 3, and the leftmost term is $i_B(\chi_{x,\vv}')^{alg}$, where $\chi_{x,\vv}'$ is the character sending 
$t = \left(\begin{smallmatrix} t_1 &  \\  & t_2 \end{smallmatrix} \right)\in T_{\GL(2)}
(\K_{\vv})$ to
$$
\chi_{x,\vv}'(t)=\chi_{x,\vv}(\begin{pmatrix} t_2 & \\ & t_1\end{pmatrix})\cdot \epsilon(\frac{t_1}{t_2}).
$$
In particular, its algebraic part $\psi_{x,\vv}'(t)=t_1^{\kappa_2}t_2^{\kappa_1-1}$ is {\it{not}} dominant, relative to the lower triangular Borel, and consequently we have $i_B(\psi_{x,\vv}')=0$. 

\begin{lem} The two extreme terms of the exact sequence are:
\begin{itemize}
\item[(a)] $i_B(\chi_{x,\vv}')^{alg}=0$,
\item[(b)] $i_B(\chi_{x,\vv})^{alg}\simeq i_B(\psi_{x,\vv})\otimes i_B(\theta_{x,\vv})$, 
\end{itemize}
where, in (b), the first induction is algebraic-induction (thus $\xi_{x,\vv}=i_B(\psi_{x,\vv})$ is irreducible algebraic of highest weight $\psi_{x,\vv}$), and the second is smooth-induction.
\end{lem}

\noindent {\it{Proof}}. For (b), note that there is a natural multiplication map,
$$
i_B(\theta_{x,\vv})\longrightarrow \Hom(\xi_{x,\vv}, i_B(\chi_{x,\vv}))^{sm},
$$
which is $\GL_2(\K_{\vv})$-equivariant. It is injective. Indeed, a function in the kernel would annihilate the highest weight vector in $\xi_{x,\vv}$, and therefore vanish on 
$\bar{N}$, and hence on the dense open Bruhat cell $B \bar{N}$. To show surjectivity, it suffices to show 
$$
i_B(\theta_{x,\vv})^K\overset{\sim}{\longrightarrow} \Hom_K(\xi_{x,\vv}, i_B(\chi_{x,\vv})),
$$
for all sufficiently small compact open subgroups $K$. All we have to do is count dimensions. On the left-hand side, we get $|B \backslash G/K|$, at least for small enough $K$.
To deal with the right-hand side, think of it as 
$$
(\xi_{x,\vv}^{\vee}\otimes i_B(\chi_{x,\vv}))^K=i_B(\xi_{x,\vv}^{\vee}|_B\otimes \chi_{x,\vv})^K.
$$
If $K$ is small enough, the latter is identified with $|B \backslash G/K|$ copies of $(\xi_{x,\vv}^{\vee}|_B\otimes \psi_{x,\vv})^B$, since $B \cap K$ is Zariski-dense, which is the line spanned by the highest weight vector in $\xi_{x,\vv}^{\vee}$. (Note here that $\xi_{x,\vv}^{\vee}$ has highest weight $\psi_{x,\vv}^{-1}$ relative to $B$.) As a result, the initial "multiplication map" is an isomorphism, and
$$
\xi_{x,\vv}\otimes i_B(\theta_{x,\vv}) \overset{\sim}{\longrightarrow} i_B(\chi_{x,\vv})^{\xi_{x,\vv}-alg}=i_B(\chi_{x,\vv})^{alg},
$$
where the last equality is deduced by suitably adapting the above argument. The same goes for (a). If $i_B(\chi_{x,\vv}')$ had $W$-algebraic vectors, for some $W$, one could deduce that $W^{\vee}$ has highest weight $\chi_{x,\vv}'^{-1}$ relative to $B$, but then $\chi_{x,\vv}'$ would be $\bar{B}$-dominant, which we already observed it is not. $\square$

\medskip

\noindent The upshot is an embedding,
$$
B(\rho_{x,\vv})^{alg} \hookrightarrow \xi_{x,\vv}\otimes i_B(\theta_{x,\vv}).
$$
On the other hand, since $\rho_{x,\vv}$ is assumed to be potentially semistable, with distinct Hodge-Tate weight, the algebraic vectors can be expressed as
$\xi_{x,\vv}\otimes \pi_{x,\vv}$, where $\pi_{x,\vv}$ arises from $WD(\rho_{x,\vv})^{F-ss}$ via the generic local Langlands correspondence (suitably normalized). We infer that
there is an inclusion of $\pi_{x,\vv}$,
$$
\pi_{x,\vv} \hookrightarrow \Hom(\xi_{x,\vv}, \xi_{x,\vv}\otimes i_B(\theta_{x,\vv}))^{sm}=(\End (\xi_{x,\vv})\otimes i_B(\theta_{x,\vv}))^{sm}=i_B(\theta_{x,\vv}),
$$
since $\xi_{x,\vv}|_K$ remains irreducible for any $K$. This proves what we want: Using Theorem 1, we have already deduced $\rho_x\simeq \rho_{\pi}$ is modular. From local-global compatibility above $p$ (due to Barnet-Lamb, Gee, Geraghty, Taylor, and Caraiani) we conclude that $\pi_{x,\vv}=BC_{\vv|v}(\pi_v)$. (Indeed, since $\rho_{\pi}$ is irreducible, $BC_{\K/F}(\pi)$ must be cuspidal, hence globally generic, and therefore $\pi_v$ is generic, so in this case "generic" local Langlands is just "classical" local Langlands).
Finally, by taking the tensor product over all $v|p$, we get the desired embedding $\pi_p \hookrightarrow i_B(\theta_x)$, which shows $x$ must be a classical point. This finishes the proof of Theorem 2, assuming Theorem 1.

\section{Proof of Theorem 1}

\subsection{Completed cohomology and Hecke algebras}

\subsubsection{Finite level Hecke algebras}

For any compact open subgroup $K \subset G(\A_f)$, the arithmetic manifold of $G$,
$$
Y(K)=G(\Q)\backslash G(\A_f)/K,
$$
is a finite set. (Recall that $G(\R)$ is compact.) For any commutative ring $A$, we let $H^0(K)_A$ be the set of functions $Y(K)  \rightarrow A$. It is naturally a module for the Hecke
algebra $\HH(K)_A$ of compactly supported $K$-biinvariant functions $G(\A_f)\rightarrow A$, equipped with convolution. We will often assume $K$ factors as a direct product 
$K=\prod_{v<\infty}K_v$, where $K_v \subset U(F_v)$ is hyperspecial for almost all places $v$ of $F$. We introduce the following finite sets of finite places of $F$, 

\begin{itemize}
\item $\Sigma_p=\{v|p\}$,
\item $\Sigma_0=\{v\nmid p: \text{$K_v$ is \underline{not} hyperspecial}\}$,
\item $\Sigma=\Sigma_0 \sqcup \Sigma_p$.
\end{itemize}

\noindent Then $\HH(K)_A$ factors as a tensor product $\otimes_{v<\infty}\HH(K_v)_A$. We will be interested in the action of the spherical part. That is, the central subalgebra below,
$$
\HH(K^p)^{sph}=\HH(K^{\Sigma})_A=\otimes_{v \notin \Sigma} \HH(K_v)_A.
$$
We let $\T(K)_A$ denote the quotient of $\HH(K^p)^{sph}$ which acts faithfully on $H^0(K)_A$. In other words, the $A$-subalgebra of $\End H^0(K)_A$ generated by all Hecke operators from $\HH(K_v)_A$, where $v \notin \Sigma$. As an $A$-module, $\T(K)_A=A \otimes \T(K)_{\Z}$ is finite free (since this holds for $A=\Z$, and hence for all rings). In particular, when $A$ is a {\it{field}}, for dimension reasons $\T(K)_A$ is therefore Artinian. (Consequently, prime ideals are maximal, there are only finitely many of them, and $\T(K)_A$ is Noetherian and semi-local; the direct product of all its localizations. When $\T(K)_A$ is reduced, these localizations can be identified with its residue fields.)

\medskip

\noindent As a Hecke-module, $H^0(K)_{\C}$ breaks up as a direct sum of (finitely many) simple modules $\pi_f^K$, where $\pi$ runs over automorphic representations of $G(\A)$,
with $\pi_{\infty}=1$. As a result, $H^0(K)_{\C}$ is a sum of simultaneous eigenspaces for $\HH(K^p)^{sph}$, and 
$\T(K)_{\C}$ is a semisimple algebra $\C \times \cdots \times \C$, where the direct factors correspond to automorphic Hecke eigensystems $\T(K)_{\C}\rightarrow \C$ (giving the action on some $\pi_f^K$). The same holds over $\bar{\Q}_p$, by transferring via an isomorphism $\C \simeq \bar{\Q}_p$.

\medskip

\noindent We will usually work over a finite extension $E/\Q_p$, with integers $\OO$, uniformizer $\varpi$, and residue field $k$. We are primarily interested in the $\OO$-algebra
$\T(K)_{\OO}$, which is known to factor as a direct product of localizations,
$$
\T(K)_{\OO}\overset{\sim}{\longrightarrow} \prod_{\m} \T(K)_{\OO,\m},
$$
with $\m$ ranging over its maximal ideals (which correspond to maximal ideals of $\T(K)_k$ via the reduction map). Here each factor $\T(K)_{\OO,\m}$ is a complete local Noetherian $\OO$-algebra. Furthermore, after extending scalars to $E$, 
$$
\T(K)_{\OO,\m} \otimes E \overset{\sim}{\longrightarrow} \prod_{\p\subset \m} \T(K)_{E,\p},
$$
where $\p\subset \m$ runs over the minimal primes of $\T(K)_{\OO}$ (which correspond to maximal ideals of $\T(K)_E$ via the inclusion map; which in turn, after fixing $E \hookrightarrow \bar{\Q}_p$, correspond to Galois-conjugacy classes of eigensystems $\T(K)_E \rightarrow \bar{\Q}_p$).

\subsubsection{Infinite level Hecke algebras}

Emerton developed the theory of completed cohomology in great generality in \cite{Em}. Here we will only use this machinery in degree zero, where everything can be done by hand. For each choice of tame level $K^p$, there are two modules,

\begin{itemize}
\item $\hat{H}^0(K^p)_{\OO}= (\underset{K_p}{\varinjlim} H^0(K_pK^p)_{\OO})^{\wedge}$ (where $\wedge$ is $p$-adic completion), 
\item $\tilde{H}^0(K^p)_{\OO}= \underset{s}{\varprojlim}  (\underset{K_p}{\varinjlim} H^0(K_pK^p)_{\OO/\varpi^s \OO}$).
\end{itemize}

\noindent There is a natural isomorphism $\hat{H}^0(K^p)_{\OO} \overset{\sim}{\longrightarrow} \tilde{H}^0(K^p)_{\OO}$, and both are naturally identified with the lattice of all continuous functions $Y(K^p) \rightarrow \OO$. Similarly we define $\hat{H}^0(K^p)_E$, which thus becomes a Banach $E$-space (for the sup-norm), endowed with a unitary 
$G(\Q_p)$-action (via right translations). Moreover, $\hat{H}^0(K^p)_E$ becomes a Banach-module for the completed Hecke algebra $\hat{\HH}(K^p)$ of biinvariant functions, which vanish at infinity (and the $G(\Q_p)$-action is Hecke-linear).

\medskip

\noindent The $p$-levels $K_p$ are ordered by reverse-inclusion; if $K_p'\subset K_p$, there is a surjective transition map
$\T(K_p'K^p)_{\OO}\twoheadrightarrow \T(K_pK^p)_{\OO}$, which makes the collection of all $\T(K_pK^p)_{\OO}$ into a projective system, as $K_p$ varies. We define
the $\OO$-algebra
$$
\T(K^p)_{\OO}=\underset{K_p}{\varprojlim} \T(K_pK^p)_{\OO},
$$
with its projective limit topology. Thus $\T(K^p)_{\OO}$ is a reduced, compact, complete $\OO$-algebra, equipped with a natural map $\HH(K^p)^{sph}\rightarrow \T(K^p)_{\OO}$ having dense image. Moreover, the action of $\HH(K^p)^{sph}$ on the completed cohomology $\hat{H}^0(K^p)_{\OO}$ extends naturally to a faithful action of $\T(K^p)_{\OO}$, as follows: Say $h=(h_{K_p})$ is a compatible sequence in $\T(K^p)_{\OO}$, and $f$ is a $p$-smooth function $Y(K^p)\rightarrow \OO$. Then $h(f)=h_{K_p}(f)$, if $f$ is $K_p$-invariant. (One verifies the right-hand side is independent of the choice of $K_p$.) This defines a continuous action on the $p$-smooth functions, which extends to the completion.
Thus,
$$
\begin{CD}
\HH(K^p)^{sph} @>dense>>\T(K^p)_{\OO} @>faithful>> End^{cts} \hat{H}^0(K^p)_{\OO}.
\end{CD}
$$
(A short argument shows that $\T(K^p)_{\OO}$ is weakly closed.)

\medskip

\noindent As in the finite-level case, $\T(K^p)_{\OO}$ is semi-local: If $K_p' \vartriangleleft K_p$ are pro-$p$ groups, an eigensystem $\HH(K^p)^{sph}\rightarrow k$ occurs in $H^0(K)_k$ if and only if it occurs in $H^0(K')_k$. Therefore $\T(K')_{\OO}\twoheadrightarrow \T(K)_{\OO}$ identifies the maximal ideals of $\T(K)_{\OO}$ with those of $\T(K')_{\OO}$.
Passing to the limit over all $K_p$, therefore yields 
$$
\text{$\T(K^p)_{\OO} \overset{\sim}{\longrightarrow} \prod_{\m} \T(K^p)_{\OO,\m}$, $\y$ $\T(K^p)_{\OO,\m}=\underset{K_p}{\varprojlim} \T(K_pK^p)_{\OO,\m}$,}
$$
with $\m$ ranging over its maximal ideals (corresponding to Galois-conjugacy classes of eigensystems occurring in $H^0(K^p)_k$). Each of the (finitely many) factors $\T(K^p)_{\OO,\m}$ is a complete local $\OO$-algebra, and they play a key role in deformation theory of Galois representations. Correspondingly,
$$
\text{$\hat{H}^0(K^p)_{\OO} \overset{\sim}{\longrightarrow} \bigoplus_{\m} \hat{H}^0(K^p)_{\OO,\m}$, $\y$ 
$\hat{H}^0(K^p)_{\OO,\m}=\hat{H}^0(K^p)_{\OO}\otimes_{\T(K^p)_{\OO}} \T(K^p)_{\OO,\m} $.}
$$
These summands $\hat{H}^0(K^p)_{\OO,\m}$ are tightly connected to the $p$-adic local Langlands correspondence, as we will see below in the course of the proof of Theorem 1.
Again, we wish to emphasize that these constructions (where $K_p$ shrinks) are due to Emerton; see Section 5.2, p. 46, in \cite{Em1}, for instance.

\subsection{Locally algebraic vectors in completed cohomology}

Inside the Banach space $\hat{H}^0(K^p)_E$, we have the dense subspace of locally analytic vectors, $\hat{H}^0(K^p)_E^{an}$, on which the Lie algebra $\g=\text{Lie} G(\Q_p)$ acts. In turn, inside the locally analytic vectors, we have the locally algebraic vectors, 
$$
\hat{H}^0(K^p)_E^{alg}=\oplus_{\xi} \hat{H}^0(K^p)_E^{\xi-alg}.
$$
Here $\xi$ runs over the absolutely irreducible algebraic $E$-representations of $G_{/E}$, up to equivalence, and the superscript $\xi-alg$ means we take the subspace of locally $\xi$-algebraic vectors; that is, those in the image of the evaluation map
$$
\xi \otimes \Hom_{K_p} (\xi,\hat{H}^0(K^p)_E) \hookrightarrow \hat{H}^0(K^p)_E,
$$
for some sufficiently small $K_p$. Equivalently, those vectors in the image of 
$$
\xi \otimes \Hom_{\g}(\xi,\hat{H}^0(K^p)_E^{an}) \hookrightarrow \hat{H}^0(K^p)_E^{an},
$$
see Proposition 4.2.10 in \cite{Em2}. Here the tensor product is over $\End_G(\xi)$, a priori a finite extension of $E$. However, in our case $G_{/E}$ is a product of copies of $\GL_2$, in which case it is easy to see $\End_G(\xi)=E$ (using the standard polynomial model of $\xi$, say).

\medskip

\noindent Each $\xi$ defines a local system $\mathcal{V}_{\xi}$ on $Y(K)$, and its space of sections $H^0(K,\mathcal{V}_{\xi})_E$ is naturally identified with the space of algebraic modular forms on $G$, of level $K$, and weight $\xi$. That is, all functions
$$
\text{$f: G(\Q)\backslash G(\A_f) \rightarrow \xi$, $\y$ $f(gk)=\xi(k_p^{-1})f(g)$, $\y$ $\forall k \in K$.}
$$
(This is straightforward; for details, see Proposition 3.2.2 in \cite{Em}, for example.) We let 
$H^0(K^p,\mathcal{V}_{\xi})_E$ denote the collection of all such $f$, where we allow $K_p$ to shrink to the identity. It carries an action of $G(\Q_p)$, and as is well-known,
$$
H^0(K^p,\mathcal{V}_{\xi})_{\C}\simeq {\bigoplus}_{\pi: \pi_{\infty}=\xi^{\vee}} m_G(\pi)\cdot \pi_p \otimes (\pi_f^p)^{K^p},
$$
where $\pi$ runs over all automorphic representations of $G(\A)$ with $\pi_{\infty}=\xi^{\vee}$. Most likely, $m_G(\pi)=1$, but we will not need that. (We tacitly move between complex coefficients and $p$-adic coefficients via a choice of isomorphism $\iota$ as above.)

\medskip

\noindent The purpose of this section is to remind the reader of the following result.

\begin{prop}
Let $\xi$ be an irreducible algebraic $E$-representation $\xi$ of $G_{/E}$. 
\begin{itemize}
\item[(a)] $H^0(K^p,\mathcal{V}_{\xi^{\vee}})_E \overset{\sim}{\longrightarrow} \Hom_{\g}(\xi,\hat{H}^0(K^p)_E^{an})$.
\item[(b)] $\xi \otimes_E H^0(K^p,\mathcal{V}_{\xi^{\vee}})_E \overset{\sim}{\longrightarrow} \hat{H}^0(K^p)_E^{\xi-alg}$.
\end{itemize}
\end{prop}

\noindent {\it{Proof}}. Clearly (b) follows from (a), in conjunction with the preceding remarks. Part (a) can really be proved by hand, so to speak (which was done in 2.2 of \cite{S}, for instance), but can more conveniently be explained as a very special case of the general machinery developed in \cite{Em}. Indeed the map in (a) is the edge map of a certain spectral sequence (given by Corollary 2.2.18 in \cite{Em}). Since $Y(K)$ is zero-dimensional, most terms vanish, and the edge map is an isomorphism. See Corollary 2.2.25 in loc. cit.,  which deduces the isomorphism in (b). $\square$ 

\begin{cor}
$\hat{H}^0(K^p)_E^{\xi-alg}\simeq {\bigoplus}_{\pi: \pi_{\infty}=\xi} m_G(\pi)\cdot (\xi\otimes \pi_p) \otimes (\pi_f^p)^{K^p}$.
\end{cor}

\noindent {\it{Proof}}. Follows from (b) and the automorphic description of $H^0(K^p,\mathcal{V}_{\xi^{\vee}})$. $\square$

\subsection{Universal modular deformations}

Suppose $\bar{\rho}: \Gamma_{\K}\rightarrow \GL_2(k)$ is an absolutely irreducible representation, which we assume to be modular of level $K=K_pK^p$, in the sense that there exists a maximal ideal $\m \subset \T(K)_k$ with which $\bar{\rho}$ is associated; by which we mean that $\m$ contains
$$
\text{$(t_w)_K-\tr \bar{\rho}(\Frob_w)$, $\y$ $(s_w)_K-\det \bar{\rho}(\Frob_w)$,}
$$
for all $w|v$ with $v \notin \Sigma$. By $(t_w)_K$ we mean the operator on $H^0(K)_k$ defined by $t_w$, and similarly for $s_w$ (see the end of Section 2.1). Note that 
$k=\T(K)_k/\m$.

\medskip

\noindent {\it{Remark}}. Eventually we will take $\bar{\rho}=\bar{\rho}_x$, for a point $x \in X(E)$ as in Theorem 1. Then for $K_p$ deep enough (pro-$p$ suffices), the reduced eigensystem $\bar{\lambda}_x$ factors through $\T(K)_k$, and we may take $\m=\m_x=\ker(\bar{\lambda}_x)$ above.

\medskip

\noindent The goal of this section is to introduce the so-called universal modular deformation of $\bar{\rho}$ (of level $K$), defined as follows.

\begin{prop}
Up to equivalence, there is a unique continuous representation
$$
\rho_{\m}=\rho_{\m,K}: \Gamma_{\K}\rightarrow \GL_2(\T(K)_{\OO,\m})
$$
such that for every place $w|v$, with $v \notin \Sigma$,
$$
\text{$\tr \rho_{\m}(\Frob_w)=(t_w)_K$, $\y$ $\det \rho_{\m}(\Frob_w)=(s_w)_K$.}
$$
Moreover,
\begin{itemize}
\item $\bar{\rho}_{\m}\simeq \bar{\rho}$,
\item For every modular deformation $\rho_{\pi,\iota}: \Gamma_{\K}\rightarrow \GL_2(A)$ of $\bar{\rho}$ (where $A$ is a complete local Noetherian $\OO$-algebra, with reside field $k$), there is a unique local morphism $\T(K)_{\OO,\m}\rightarrow A$ with respect to which $\rho_{\pi,\iota}$ is the specialization of $\rho_{\m}$.
\end{itemize}
(In this sense, $\rho_{\m}$ is the universal level-$K$ modular deformation of $\bar{\rho}$.)
\end{prop}

\noindent {\it{Proof}}. First, for each minimal prime $\p \subset \m$ of $\T(K)_{\OO}$, we define a representation
$$
\text{$\rho_{\p}: \Gamma_{\K}\rightarrow \GL_2(\T(K)_{E,\p})$, $\y$ $\T(K)_{E,\p}=\T(K)_{E}/\p$,}
$$ 
as follows: Pick an eigensystem $\lambda:\T(K)_E \rightarrow \bar{\Q}_p$ with $\p=\ker(\lambda)$, extending our chosen embedding $E \hookrightarrow \bar{\Q}_p$. Its restriction to $\T(K)_{\Q}$ then arises from an automorphic representation $\pi$, of level $K$, and weight $\pi_{\infty}=1$. That is, $\lambda=\iota\lambda_{\pi}$. After possibly enlarging $E$, 
we may assume $\rho_{\pi,\iota}$ takes values in $\GL_2(\OO)$. We then let $\rho_{\p}=\rho_{\pi,\iota}$, and observe (by Chebotarev) that this only depends on $\p$; not the choice of $\pi$. By construction, for all $w|v$ with $v \notin \Sigma$,
$$
\text{$\tr \rho_{\p}(\Frob_w)=(t_w)_K+\p$, $\y$ $\det \rho_{\p}(\Frob_w)=(s_w)_K+\p$.}
$$
By the factorization at the end of 4.1.1, we obtain a representation
$$
\varrho=\prod_{\p\subset \m}\rho_{\p}: \Gamma_{\K}\rightarrow \GL_2(\T(K)_{\OO,\m}\otimes E)
$$
with similar properties. In particular, $\tr\varrho(\Frob_w)=(t_w)_K$, which shows that the pseudo-character $\tr \varrho$ takes values in the complete local ring $\T(K)_{\OO,\m}$.
Its reduction equals $\tr \bar{\rho}$; the trace of an absolutely irreducible representation. By a result of Nyssen, Rouquier, and Procesi, there is a unique representation
$$
\rho_{\m}=\rho_{\m,K}: \Gamma_{\K}\rightarrow \GL_2(\T(K)_{\OO,\m})
$$
with trace $\tr \varrho$. (We refer to Chapter 1 of \cite{BC} for an in-depth treatment of pseudo-characters and for when they originate from actual representations.) $\square$

\medskip

\noindent Finally, we will let $K_p$ shrink. When $K_p$ varies, the representations $\rho_{\m,K_pK^p}$ are compatible, and we may pass to the limit, resulting in
$$
\rho_{\m}=\rho_{\m,K^p}=\underset{K_p}{\varprojlim}  \rho_{\m,K_pK^p}:\Gamma_{\K}\rightarrow \GL_2(\T(K^p)_{\OO,\m}),
$$
the universal modular deformation of $\bar{\rho}$, of tame level $K^p$. Below, we shall link its local restrictions $\rho_{\m,\vv}$ to the $p$-adic local Langlands correspondence
for $\GL_2(\Q_p)$, using the deformation-theoretic approach of Kisin \cite{Ki}.

\subsection{Colmez's Montreal functor}

It turns out it is easier to describe the inverse of the $p$-adic local Langlands correspondence $B(\cdot)$. At a 2005 workshop in Montreal, Colmez gave an elegant definition of a functor $\VV$, from representations of $\GL_2(\Q_p)$ to those of $\Gamma_{\Q_p}$, which he subsequently studied in detail \cite{Col}. More precisely, $\VV$ is an exact covariant functor,
$$
\VV: \Rep_{\OO}(\GL_2(\Q_p)) \longrightarrow \Rep_{\OO}(\Gamma_{\Q_p}),
$$
from the category of smooth, admissible, finite length representations of $\GL_2(\Q_p)$ on torsion $\OO$-modules (which admit a central character) to the category 
of representations of $\Gamma_{\Q_p}$ on finite length $\OO$-modules. 

\begin{thm}
The functor $\VV$ enjoys the following properties.
\begin{itemize}
\item[(1)] $\VV(\pi)=0$ if and only if $\pi$ is finite over $\OO$.
\item[(2)] $\VV$ is compatible with twisting by a character.
\item[(3)] If $V$ is a two-dimensional $k$-representation of $\Gamma_{\Q_p}$, which is \underline{not} a twist of 
$$
\begin{pmatrix}1 & * \\ & \bar{\epsilon}\end{pmatrix},
$$
then there exists a unique representation $\bar{\pi}\in \Rep_{k}(\GL_2(\Q_p))$ such that
\begin{itemize}
\item[(a)] $\VV(\bar{\pi})\simeq V$, and $\bar{\pi}$ has central character (corresponding to) $\det(V)\bar{\epsilon}^{-1}$.
\item[(b)] $\Ext_{\GL_2(\Q_p)}^1(\bar{\pi},\bar{\pi})\hookrightarrow \Ext_{\Gamma_{\Q_p}}^1(V,V)$.
\item[(c)] $\bar{\pi}^{\SL_2(\Q_p)}=0$.
\end{itemize}
(Thus $\bar{\pi}$ has no finite-dimensional submodules or quotients.)
\end{itemize} 
Furthermore, $\VV$ realizes the mod $p$ local Langlands correspondence: $\bar{\pi}^{ss}\leftrightarrow V^{ss}$.
\end{thm}

\noindent {\it{Proof}}. This is taken from \cite{Ki}, where it occurs as Theorem 2.1.1. It summarizes some of the most important results from \cite{Col}. $\square$

\medskip

\noindent The representation $\bar{\pi}$ can be written down explicitly when $V$ is reducible; see Remark 3.3.3 in \cite{Em1}. For instance, if $V$ is a sum of characters
(in general position), $\bar{\pi}$ is a sum of two irreducible principal series. In general, $\bar{\pi}$ sits in an extension of such. On the other hand, when $V$ is irreducible, $\bar{\pi}$ is supersingular.

\medskip

\noindent  Although we will not use the actual definition of $\VV$, we find it worthwhile to briefly outline its construction. We will follow Section 3.1 in \cite{Bre}. Thus, we let $\pi$ be a smooth, admissible, finite length representations of $\GL_2(\Q_p)$ on an $\OO/\varpi^s\OO$-module. Choose a $\GL_2(\Z_p)$-stable $\OO$-submodule $W \subset \pi$ such that 
$$
\pi=\OO[\GL_2(\Q_p)]\cdot W.
$$
(Such $W$ exist, since $\pi$ has finite length; use induction on the length.) Introduce
$$
W^{\#}=\sum_{m=0}^{\infty}\begin{pmatrix}p^m & \Z_p \\ & 1\end{pmatrix}\cdot W,
$$
the smallest $\OO$-submodule $W^{\#}\subset \pi$, containing $W$,  which is stable under the actions of $\left(\begin{smallmatrix} p & \\ & 1\end{smallmatrix} \right)$ and 
$\left(\begin{smallmatrix} 1 & \Z_p\\ & 1\end{smallmatrix}\right)$; the latter acting smoothly. In addition, $\left(\begin{smallmatrix} \Z_p^{\times} & \\ & 1\end{smallmatrix}\right)$ preserves $W^{\#}$. Next, we look at its Pontryagin dual. That is, 
$$
M(W)=\Hom_{\OO/\varpi^s\OO}(W^{\#},\OO/\varpi^s\OO),
$$
which then naturally becomes a module over the completed group ring,
$$
\text{$\OO[[\begin{pmatrix} 1 & \Z_p\\ & 1\end{pmatrix}]]\simeq \OO[[X]]$, $\y$ $\begin{pmatrix}1 & 1 \\ & 1\end{pmatrix} \leftrightarrow X+1$.}
$$
Moreover, $M(W)$ carries an action of $\begin{pmatrix} p & \\ & 1\end{pmatrix}$ and $\Gamma\simeq \Z_p^{\times}$. Consequently,
$$
D=M(W) \otimes_{\OO[[X]]} \OO((X))
$$
comes with the structure of a $(\phi,\Gamma)$-module: A finite type $\OO((X))$-module, annihilated by $\varpi^s$, endowed with commuting semilinear actions of $\phi$ and $\Gamma$, such that $\phi(D)$ generates $D$. In \cite{Col}, Colmez verifies that this $D$ is independent of the choice of $W$. Finally, by a fundamental result of Fontaine, the category of such $(\phi,\Gamma)$-modules is equivalent to the category of representations of $\Gamma_{\Q_p}$ on finite-length $\OO/\varpi^s\OO$-modules. One defines 
$\VV(\pi)^{\vee}$ to be the representation corresponding to $D$; we take the contragredient in order to make $\VV$ covariant. 

\medskip

\noindent {\it{Remark}}. Note that $\VV(\pi)$ need not be two-dimensional, if $\pi$ is a $k$-representation of $\GL_2(\Q_p)$. For example, if $\pi$ is a principal series (or twisted Steinberg), $\VV(\pi)$ is a character. We refer to the bottom of p. 19 in \cite{Em1} for a table of values of $\VV$.

\medskip

\noindent Finally, let us point out that $\VV$ can be defined for more general coefficient rings $A$. The above recipe works, more or less ad verbatim, for a local Artinian $\OO$-algebra $A$ whose residue field is a finite extension of $k$. More generally, if $A$ is a complete local Noetherian $\OO$-algebra (with residue field finite over $k$), and $\pi$ is a "suitable" $A$-representation of $\GL_2(\Q_p)$, one defines $\VV(\pi)$ as the inverse limit of $\VV(\pi/\m_A^s \pi)$. See Section 3.2 of \cite{Em1} for a more thorough discussion of this.

\subsection{Deformation theory and $p$-adic Langlands}

We go back to our absolutely irreducible, modular, representation $\bar{\rho}$ from 4.2. At each place $w|p$, we look at its restriction $\bar{\rho}_{w}$ to the Galois group of $\K_w \simeq \Q_p$. One of our standing hypotheses is that (for all $\chi$ and extensions $*$),
$$
\bar{\rho}_{w} \nsim \chi \otimes \begin{pmatrix}1 & * \\ & \bar{\epsilon}\end{pmatrix}.
$$
Therefore, $\bar{\rho}_{w}\leftrightarrow \bar{\pi}_{w}$, a smooth, admissible, finite length $k$-representation of $\GL_2(\Q_p)$, via mod $p$ local Langlands. Since $\VV$ is 
exact, it takes a deformation of $\bar{\pi}_{w}$ to a deformation of $\bar{\rho}_{w}$. Recall that a deformation of $\bar{\rho}_w$ to a complete local Noetherian $\OO$-algebra $A$
(with residue field $k$) is a free rank two $A$-module $\rho_w$, with a continuous $\Gamma_{\K_w}$-action, such that $\rho_w \otimes_A k\simeq \bar{\rho}_w$. Deformations of $\bar{\pi}_w$ are defined similarly; we refer to Definition 3.2 in \cite{Bre}. 

\medskip

\noindent For simplicity only, we will assume $\End_{\Gamma_{\K_w}}(\bar{\rho}_w)=k$ (for example, this holds 
if $\bar{\rho}_w$ is absolutely irreducible). Equivalently, $\End_{\GL_2(\K_w)}(\bar{\pi}_w)=k$. This forces a certain rigidity into our deformation problems, and Schlessinger's criterion guarantees they are representable, by complete local Noetherian $\OO$-algebras $R(\bar{\rho}_w)$ and $R(\bar{\pi}_w)$, respectively. One can relax the "endomorphism-condition" and work with a more advanced deformation theory (making use of groupoid-valued functors, and framings). We refrain from doing so. The necessary modifications go exactly as on p. 25 in \cite{Em1}.

\medskip

\noindent Since $\VV$ is exact, it defines a morphism of local $\OO$-algebras,
$$
\VV: R(\bar{\rho}_w)\longrightarrow R(\bar{\pi}_w),
$$
which we will continue to denote $\VV$. Let $R(\bar{\pi}_w)^{\det}$ denote the quotient of $R(\bar{\pi}_w)$, which parametrizes deformations $\pi_w$, which admit a central character $\chi_{\pi_w}$ corresponding to $\det \VV(\pi_w)\epsilon$, via local class field theory. We say $\pi_w$ satisfies the "determinant-condition". In \cite{Col}, VII.5.3, Colmez showed the composition
$$
R(\bar{\rho}_w)\overset{\VV}{\longrightarrow} R(\bar{\pi}_w) \longrightarrow R(\bar{\pi}_w)^{\det}
$$
is {\it{onto}}. Geometrically, $\Spec R(\bar{\pi}_w)^{\det}$ is a closed subset of the deformation space $\Spec R(\bar{\rho}_w)$. We will intersect it with another closed subset; the Zariski-closure of the crystalline points. More precisely, consider the quotient
$$
R(\bar{\rho}_w)^{cris}=R(\bar{\rho}_w)/\cap\p,
$$
where $\p$ runs over all prime ideals of the form $\p=\ker(R(\bar{\rho}_w)\overset{\alpha}{\rightarrow} \bar{\Z}_p)$, where $\alpha$ is a homomorphism such that the $\alpha$-specialization is a crystalline representation $\Gamma_{\K_w}\rightarrow \GL_2(\bar{\Q}_p)$ with distinct Hodge-Tate weights. We say $\p$ runs over the "crystalline-regular" points of $\Spec R(\bar{\rho}_w)$. Thus we may think of $\Spec R(\bar{\rho}_w)^{cris}$ as their Zariski-closure. We look at the "intersection",
$$
\Spec R(\bar{\pi}_w)^{cris}=\Spec R(\bar{\rho}_w)^{cris} \times_{\Spec R(\bar{\rho}_w)} \Spec R(\bar{\pi}_w)^{\det}.
$$
In other words, the tensor product $R(\bar{\pi}_w)^{cris}$ fits in a Cartesian square,
$$
\begin{CD}
R(\bar{\rho}_w) @>>> R(\bar{\pi}_w)^{\det}\\
@VVV      @VVV\\
R(\bar{\rho}_w)^{cris}  @>>> R(\bar{\pi}_w)^{cris}
\end{CD}
$$
Here the bottom map turns out to be an isomorphism, by a key result of Kisin.

\begin{thm}
$R(\bar{\rho}_w)^{cris} \overset{\sim}{\longrightarrow} R(\bar{\pi}_w)^{cris}$.
\end{thm}

\noindent {\it{Proof}}. This is Proposition 2.3.3 in \cite{Ki}, and its Corollary 2.3.4. $\square$

\medskip

\noindent Intuitively, this says that $\Spec R(\bar{\rho}_w)^{cris}$ is contained in $\Spec R(\bar{\pi}_w)^{\det}$. Hence, what goes into the proof is to first show that crystalline-regular deformations $\rho_w$ lie in the image of $\VV$, say $\rho_w=\VV(\pi_w)$; and thereafter that such $\pi_w$ automatically satisfies the determinant-condition.

\medskip

\noindent We will now apply this to our universal modular deformation $\rho_{\m}$ over $\T(K^p)_{\OO,\m}$. Its restriction $\rho_{\m,w}$ is a deformation of $\bar{\rho}_w$, so it arises from the universal deformation over $R(\bar{\rho}_w)$ via a unique morphism of local $\OO$-algebras,
$$
R(\bar{\rho}_w) \overset{\alpha}{\longrightarrow} \T(K^p)_{\OO,\m}.
$$ 
We wish to show $\alpha$ factors through $R(\bar{\rho}_w)^{cris}$. This, combined with the previous Theorem, would show the existence of a unique deformation $\pi_{\m,w}$ of
$\bar{\pi}_w$ over $\T(K^p)_{\OO,\m}$, which satisfies the determinant-condition, and such that 
$$
\rho_{\m,w}=\VV(\pi_{\m,w}).
$$
Here we will show that $\alpha$ factors, granted that the crystalline points are Zariski dense; which is the main result of the next section.

\begin{lem}
The morphism $R(\bar{\rho}_w) \overset{\alpha}{\longrightarrow} \T(K^p)_{\OO,\m}$ factors through $R(\bar{\rho}_w)^{cris}$.
\end{lem}

\noindent {\it{Proof}}. Suppose $r \in \cap \p$. We must show $\alpha(r)$ acts trivially on $\hat{H}^0(K^p)_{\OO,\m}$. By the main Proposition of the next section, 
it suffices to show $\alpha(r)$ acts trivially on each $\lambda$-eigenspace of $\hat{H}^0(K^p)_{\OO,\m}\otimes E$, with $\lambda=\lambda_{\pi}$ as in that Proposition.
For the remainder of this proof, let $\q=\ker(\lambda)$, viewed as a prime ideal in $\T(K^p)_{\OO,\m}$. We need to show its pullback $\alpha^{-1}(\q)$ is among the $\p$'s in the intersection $\cap \p$. However, the $\q$-specialization $\rho_{\m}(\q)$ of the universal modular deformation, can be identified with $\rho_{\pi,\iota}$, once we fix an embedding 
$\T(K^p)_{\OO}/\q \hookrightarrow \bar{\Q}_p$. Since $\rho_{\pi,\iota}$ is a deformation of $\bar{\rho}$, which is crystalline-regular at every $w|p$, the morphism
$$
R(\bar{\rho}_w) \rightarrow \T(K^p)_{\OO,\m} \rightarrow \T(K^p)_{\OO}/\q \hookrightarrow \bar{\Q}_p
$$ 
factors through $R(\bar{\rho}_w)^{cris}$, which is to say that $\cap \p \subset \alpha^{-1}(\q)$, as wanted. $\square$

\medskip

\noindent {\it{Remark}}. In fact, one can show that the projection $R(\bar{\rho}_w)\rightarrow R(\bar{\rho}_w)^{cris}$ is an isomorphism (that is, $\cap \p=0$ in the above notation). 
For instance, see Theorem A in \cite{Ch}, when $\bar{\rho}_w$ is irreducible. Consequently, so is $R(\bar{\pi}_w)^{det}\rightarrow R(\bar{\pi}_w)^{cris}$, and 
$$
\VV: R(\bar{\rho}_w)\overset{\sim}{\longrightarrow} R(\bar{\pi}_w)^{det}
$$
is a strengthening of Theorem 7. We note that an alternative approach to some of these results can be found in the important recent work of Paskunas \cite{Pas}.

\subsection{Zariski density of crystalline points}

First, let us get the definition of "classical and crystalline" in place.

\begin{df}
A prime ideal $\p \subset \T(K^p)_{\OO}$ is called \underline{classical} if $\p=\ker(\lambda)$, for some eigensystem $\lambda:  \T(K^p)_{\OO} \rightarrow \bar{\Q}_p$ associated with an automorphic representation $\pi$ of $G(\A)$, of tame level $K^p$ (and possibly non-trivial weight). If moreover $\pi_p$ is unramified, we say $\p$ is classical and crystalline. We denote by $\CC$ the set of all classical and crystalline points $\p$ in $\Spec \T(K^p)_{\OO}$.
\end{df}

\noindent Note that prime ideals $\p$ in the localization $\T(K^p)_{\OO,\m}$ correspond to prime ideals $\p \subset \T(K^p)_{\OO}$ such that $\p \subset \m$. We will pass between these points of view, with no mention. We will let $\CC_{\m}$ denote the set of $\p \in \CC$ contained in $\m$.

\medskip

\noindent Our main result in this section is:

\begin{prop}
The submodule $\oplus_{\p \in \CC} \hat{H}^0(K^p)_E[\p]^{alg}$ is dense in $\hat{H}^0(K^p)_E$. (Similarly, $H=\hat{H}^0(K^p)_{\OO,\m}\otimes E$ contains 
$\oplus_{\p \in \CC_{\m}} H[\p]^{alg}$ as a dense submodule.)
\end{prop}

\noindent {\it{Proof}}. This was proved in great detail in \cite{Sor}; see Corollary 4 in section 7.5. Here we will only be cursory, and content ourselves with describing the main threads of the proof, which is a variation of ideas of Katz. To simplify the argument even further, we will only treat $\CC$, and leave it to the reader to adapt this to $\CC_{\m}$.

\medskip

\noindent Throughout the proof, we will let $A=\OO/\varpi^s \OO$, for a fixed $s>0$. One first observes that $H^0(K^p,A)$ is an injective admissible $A[K_p]$-module (for sufficiently small $K^p$, and arbitrary $K_p$). That is, by definition, that the functor
$$
M \mapsto \Hom_{A[K_p]}(M, H^0(K^p,A))
$$
is exact. Indeed, the right-hand side can be interpreted as $M^{\vee}$-valued modular forms, and hence be identified with $(M^{\vee})^{\oplus h}$, for a suitable $h$.
(Here $M^{\vee}$ denotes the Pontryagin dual, which naturally becomes a module over the Iwasawa algebra $A[[K_p]]$, when $M$ is smooth.) Consequently, 
$H^0(K^p,A)^{\vee}$ is a projective finitely generated $A[[K_p]]$-module.

\medskip

\noindent Now, if we assume $K_p$ is a pro-$p$ subgroup, a more or less standard argument shows that $A[[K_p]]$ is a (non-commutative) local ring; see 7.3 in \cite{Sor} for all the gory details. As is easily checked, Nakayama's lemma still holds for non-commutative rings $R$, and in particular one infers that projective finitely generated $R$-modules are free. Applying this to our setup, with $R=A[[K_p]]$, for a pro-$p$ group $K_p$,
$$
H^0(K^p,A)^{\vee} \simeq A[[K_p]]^{r_s}
$$
as $A[[K_p]]$-modules, for some $r_s>0$, which will be shown to be independent of $s$: Taking the Pontryagin dual on both sides of the previous isomorphism, 
$$
H^0(K^p,A)\simeq \CC(K_p,A)^{r_s}.
$$
Scale both sides by $\varpi$, and compare quotients. Then take $K_p'$-invariants, for some $K_p' \vartriangleleft K_p$, and compare dimensions over $k$. On the left, one gets 
the cardinality of the set $Y(K_p'K^p)$, which is independent of $s$; and on the right one gets  $[K_p:K_p']r_s$. We deduce that $r_s=r_1$ is independent of $s$, and we simply denote it by $r$. Going back to the previous isomorphism, we can now take the inverse limit over $s$ on both sides, which results in an isomorphism
$$
\tilde{H}^0(K^p)_{\OO} \simeq \CC(K_p,\OO)^r
$$
of $\OO[K_p]$-modules (recall that $K_p$ is assumed to be pro-$p$). We remind the reader that $\tilde{H}^0(K^p)_{\OO}$ is canonically isomorphic to
$\hat{H}^0(K^p)_{\OO}$, as mentioned in 4.1.2.

\medskip

\noindent We will use the last isomorphism to show the locally $G(\Z_p)$-algebraic vectors are dense in $\hat{H}^0(K^p)_E$. First off, we pick a pro-$p$ subgroup $K_p \vartriangleleft G(\Z_p)$. Thus,
$$
\hat{H}^0(K^p)_E^{\vee}\simeq E[[K_p]]^r 
$$
is a finite free $E[[K_p]]$-module (by what we have just shown). Consequently,
$$
\Hom_{E[[G(\Z_p)]]}(\hat{H}^0(K^p)_E^{\vee},-)=\Hom_{E[[K_p]]}(\hat{H}^0(K^p)_E^{\vee},-)^{G(\Z_p)/K_p}
$$
is an exact functor, which is to say $\hat{H}^0(K^p)_E^{\vee}$ is a projective $E[[G(\Z_p)]]$-module., and therefore a direct summand of a finite free module: For some module $Z$,
$$
\hat{H}^0(K^p)_E^{\vee}\oplus Z \simeq E[[G(\Z_p)]]^e.
$$
Undoing the dual, we find $\hat{H}^0(K^p)_E$ sitting as a direct summand of $\CC(G(\Z_p),E)^e$. Observe that $G(\Z_p)$ is identified with $\prod_{v|p}\GL_2(\OO_{\vv})$, which (topologically) is a closed-open subset of $\Z_p^{4d}$. Hence, any continuous function on $G(\Z_p)$ extends (non-uniquely) to $\Z_p^{4d}$, and this extension has a multi-variable Mahler expansion, which shows it can be approximated uniformly by polynomial functions on $\Z_p^{4d}$. The latter obviously restrict to $G(\Z_p)$-algebraic functions in
$\CC(G(\Z_p),E)$.

\medskip

\noindent To finish the proof, notice that (as in the proof of Corollary 1),
$$
\hat{H}^0(K^p)_E^{G(\Z_p)-alg}\simeq \oplus_{\xi}\oplus_{\pi: \pi_{\infty}=\xi} m_G(\pi) \cdot (\xi\otimes \pi_p^{G(\Z_p)})\otimes (\pi_f^p)^{K^p},
$$
 which we now know is dense in $\hat{H}^0(K^p)_E$. A fortiori, so is the $G(\Q_p)$-submodule
 $$
\langle \hat{H}^0(K^p)_E^{G(\Z_p)-alg} \rangle \simeq \oplus_{\xi}\oplus_{\pi: \pi_{\infty}=\xi, \pi_p^{G(\Z_p)}\neq 0} m_G(\pi) \cdot (\xi\otimes \pi_p)\otimes (\pi_f^p)^{K^p},
$$
 it generates. Visibly, the only eigensystems contributing to this sum are those $\lambda: \T(K^p)_{\OO}\rightarrow E$ associated with automorphic $\pi$, of tame level $K^p$, some weight $\xi$, such that $\pi_p$ is unramified. In other words, such that $\p=\ker(\lambda)\in \CC$. We finish off by noting that the right-hand side therefore clearly lies in 
 the submodule $\oplus_{\p \in \CC} \hat{H}^0(K^p)_E[\p]^{alg}$ of the Proposition, which is therefore dense. $\square$
 
 \medskip
 
\noindent The following immediate consequence will be quite useful (in 4.9):
 
\begin{cor}
$\cap_{\p \in \CC_{\m}}\p=0$ (that is, $\CC_{\m}$ is Zariski dense in $\Spec \T(K^p)_{\OO,\m}$).
\end{cor}

\noindent {\it{Proof}}. Any $t \in \cap_{\p \in \CC_{\m}}\p$ must act trivially on all of $\hat{H}^0(K^p)_{\OO,\m}$, so $t=0$. $\square$ 

\medskip

\noindent {\it{Remark}}. The Zariski-density of crystalline points have been proved by one of us (P.C.) in the generality of a PEL-type Shimura variety, see \cite{Cho}. The argument is much similar to the zero-dimensional case, except for a few issues with cohomology in higher degree. We decided to include the above sketch (of the proof of Proposition 3) to make our text more self-contained. 

\subsection{Reformulation of Theorem 1}

We are finally in possession of all the ingredients needed to define a certain module $\XX$, in terms of which Theorem 1 gets a simple formulation. For each $w|p$, we have $\pi_{\m,w}$, a deformation of $\bar{\pi}_w \leftrightarrow \bar{\rho}_w$ over $\T(K^p)_{\OO,\m}$, with central character $\det(\rho_{\m,w})\epsilon$, such that $\rho_{\m,w}=\VV(\pi_{\m,w})$. First, we introduce
$$
\pi_{\m}=\otimes_{v|p}\pi_{\m,\vv}.
$$
(Here $v|p$ varies over places of $F$, and $\vv|v$ is our choice of a place of $\K$ above $v$. The tensor product is over the ring $\T(K^p)_{\OO,\m}$.) This $\pi_{\m}$ is thus a
$\T(K^p)_{\OO,\m}$-module, with a linear action of the group $G(\Q_p)$, which we always identify with $\prod_{v|p}\GL_2(\K_{\vv})$, using our collection of places $\{\vv\}$.
Furthermore, $\pi_{\m}$ has a natural $\m$-adic topology, from $\T(K^p)_{\OO,\m}$. Strongly inspired by Section 6.3 in \cite{Em1}, we put
$$
\XX=\XX_{K^p}=\Hom_{\T(K^p)_{\OO,\m}[G(\Q_p)]}^{cts.}(\pi_{\m}, \hat{H}^0(K^p)_{\OO,\m}).
$$
(See also Section 4.1 in \cite{Bre}). Here the big difference with \cite{Em1} and \cite{Bre}, is the lack of a Galois-action on $\hat{H}^0(K^p)$, so we look at continuous homomorphism out of $\pi_{\m}$, as opposed to $\rho_{\m}\otimes \pi_{\m}$. Moreover, we find it simpler to work with a fixed tame level $K^p$ throughout (and hence a fixed 
eigenvariety).

\medskip

\noindent We find it useful to spell out the continuity-assumption: Thus, the module $\XX$ consists of $\T(K^p)_{\OO,\m}$-linear, $G(\Q_p)$-equivariant, homomorphisms of the form
$$
\eta: \pi_{\m}\rightarrow \hat{H}^0(K^p)_{\OO,\m},
$$
such that $\forall s \in \Z_{>0}$, there is a $t \in \Z_{>0}$, such that
$$
\eta(\m^t \pi_{\m})\subset \varpi^s \hat{H}^0(K^p)_{\OO,\m}.
$$
(The reader may want to compare this to the first paragraph of 4.4 in [Bre].) 

\begin{thm}
$\XX[\p]\neq 0$, for \underline{all} prime ideals $\p\subset \T(K^p)_{\OO,\m}$.
\end{thm}

\noindent This implies our main result.

\begin{lem}
Theorem 8 $\Rightarrow$ Theorem 1.
\end{lem}

\noindent {\it{Proof}}. If $\p$ is a prime, the $\p$-torsion $\XX[\p]$ consists of those $\eta$ 
which factor through 
$$
\pi_{\m}(\p)=\pi_{\m}/\p\pi_{\m}\simeq \otimes_{v|p} (\pi_{\m,\vv}/\p \pi_{\m,\vv}),
$$
this tensor product being over the field $\T(K^p)_{\OO,\m}/\p$. Since $\VV$ is an exact functor,
$$
\rho_{\m,\vv}(\p)=\rho_{\m,\vv}/\p\rho_{\m,\vv}=\VV(\pi_{\m,\vv}/\p \pi_{\m,\vv}).
$$
If we take $\p=\p_x$, as in Theorem 1, the left-hand side is $\rho_{x,\vv}$. Consequently,
$$
\pi_{\m}(\p_x)=\otimes_{v|p} B(\rho_{x,\vv}),
$$
since $B$ and $\VV$ are each others inverse (see Definition 3.3.15 on p. 26 in \cite{Em1}). In conclusion, the non-vanishing $\XX[\p_x]\neq 0$ of Theorem 8 amounts to the existence of a nonzero, continuous, $E$-linear, $G(\Q_p)$-equivariant, homomorphism
$$
\eta: \pi_{\m}(\p_x)=\otimes_{v|p} B(\rho_{x,\vv}) \longrightarrow \hat{H}^0(K^p)_{\OO,\m}[\p_x]\simeq \hat{H}^0(K^p)_E[\p_x]
$$
(Since $\pi_{\m}(\p_x)$ is annihilated by $\p_x$, so is the image of $\eta$.) Finally, since the target is complete, and $\eta$ is continuous, it extends uniquely to the completion
(with respect to the tensor product norm). This proves Theorem 1. $\square$

\subsection{Non-vanishing at classical crystalline points}

Recall that $\bar{\rho}: \Gamma_{\K}\rightarrow \GL_2(k)$ is an absolutely irreducible representation, associated with the maximal ideal $\m \subset \T(K^p)_{\OO}$. In this section we will make the rather bold assumption that all its restrictions $\bar{\rho}_w: \Gamma_{\K_w}\rightarrow \GL_2(k)$ remain absolutely irreducible, for $w|p$. This is to avoid having to deal with $B(\rho)$, for reducible crystalline $\rho$. We are hopeful that one can adapt the approach of \cite{BE}, and relax this condition. For now, we wish to keep things simple.

\medskip

\noindent In the next section, by a formal Nakayama-type argument, we will verify that it is enough to prove Theorem 8 for a Zariski dense subset of $\p\subset \m$. Here we will prove the non-vanishing for all classical and crystalline points $\p$.

\begin{lem}
$\XX[\p]\neq 0$, for all $\p\in \CC_{\m}$.
\end{lem}

\noindent {\it{Proof}}. Let $\p \subset \m$ be a classical crystalline point of $\Spec \T(K^p)_{\OO}$. Say, $\p=\ker(\lambda_{\pi})$, for an automorphic $\pi$, which is unramified at $p$ (and of tame level $K^p$, and some weight). By the observations made in the proof of Lemma 3, we seek a nonzero, continuous, $E$-linear, $G(\Q_p)$-equivariant map,
$$
\eta: \otimes_{v|p} B(\rho_{\pi,\vv})\simeq \pi_{\m}/\p\pi_{\m}\longrightarrow \hat{H}^0(K^p)_{\OO,\m}[\p]\simeq \hat{H}^0(K^p)_E[\p],
$$
after possibly passing to a finite extension of $E$ (containing the residue field of $\p$). Here $\rho_{\pi}=\rho_{\pi,\iota}$ is as in Theorem 3; we suppress $\iota$. Note that 
$\bar{\rho}_{\pi}\simeq \bar{\rho}$, since $\p \subset \m$. In particular, by our "bold" assumption on $\bar{\rho}_{\vv}$, we infer that the crystalline representation $\rho_{\pi,\vv}$ is absolutely irreducible
(in fact, residually). 

\medskip

\noindent Now, for absolutely irreducible crystalline $\rho$, the $p$-adic local Langlands correspondent $B(\rho)$ has a simple description, due to Berger and Breuil,
\cite{BB}, which we briefly recall: Following the recipe of \cite{BS}, one first associates a locally algebraic representation $LA(\rho)$ of $\GL_2(\Q_p)$ out of the $p$-adic Hodge theoretical data of $\rho$.
$$
LA(\rho)=\xi(\rho)\otimes_E \pi(\rho)={\det}^{\kappa_1} \text{Sym}^{\kappa_2-\kappa_1-1}(E^2)\otimes_E \pi(\rho),
$$
where $\kappa_1<\kappa_2$ are the Hodge-Tate weights of $\rho$, and the smooth factor $\pi(\rho)$ is given by the generic local Langlands correspondence.
Thus, $\pi(\rho)$ is a full unramified principal series, possibly reducible. When it is irreducible, one has
$$
WD(\rho)^{F-ss}\simeq rec(\pi(\rho)\otimes |\det|^{-1/2}).
$$
In the reducible case, $WD(\rho)$ corresponds to the Langlands quotient of $\pi(\rho)$. By Theorem 2.3.2 in \cite{Ber}, which summarizes some of the main results of \cite{BB}, one knows that $LA(\rho)$ admits a separated $\GL_2(\Q_p)$-stable $\OO$-lattice, which is {\it{finitely generated}} over $\GL_2(\Q_p)$. Clearly all such lattices are commensurable, and $B(\rho)$ is the completion of $LA(\rho)$ with respect to any one of them. Thus $B(\rho)$ becomes a topologically irreducible, admissible, unitary Banach $E$-space representation.

\medskip

\noindent Specializing this discussion to $\rho=\rho_{\pi,\vv}$, Theorem 3 (especially local-global compatibility at the places above $p$) allows us to compute $LA(\rho)$ in terms of $\pi$. 
$$
\text{$\xi(\rho)=\xi_{\vv}$, $\y$ $\pi(\rho)=BC_{\vv|v}(\pi_v)$,}
$$
where $\xi_{\vv}$ denotes the irreducible algebraic representation of $\GL_2(\K_{\vv})$, over $E$, related to $\pi_{\infty}$ as follows: A priori, $\pi_{\infty}$ is a representation of $G(\C)$, restricted to $G(\R)$. Via $\iota$, we view $\pi_{\infty}$ as a representation of $G(\bar{\Q}_p)$, and restrict it to $G(\Q_p)$. The resulting representation is 
$\xi=\otimes_{v|p}\xi_{\vv}$. For more details, see 2.4 in \cite{Sor}. So,
$$
\otimes_{v|p}LA(\rho_{\pi,\vv})\simeq \xi \otimes \pi_p,
$$
both viewed as representations of $G(\Q_p)$, identified with $\prod_{v|p}\GL_2(\K_{\vv})$. 

\medskip

\noindent We now invoke Corollary 1 from Section 4.2, which shows the existence of embeddings,
$$
\otimes_{v|p}LA(\rho_{\pi,\vv})\simeq \xi \otimes \pi_p \hookrightarrow \hat{H}^0(K^p)_E[\p]^{alg},
$$
parametrized by $(\pi_f^p)^{K^p}$, and the multiplicity of $\pi$, if greater than one. 

\medskip

\noindent If $\LL$ is an arbitrary Banach $E$-space, with a unitary action of $GL_2(\K_{\vv})$, then any equivariant map 
$i:LA(\rho_{\pi,\vv})\rightarrow \LL$ is automatically continuous; with respect to the topology given by a finite-type lattice $\Lambda$ in $LA(\rho_{\pi,\vv})$. This is almost immediate. If $\Lambda$ is generated by $\{\lambda_j\}$ as a $\GL_2(\K_{\vv})$-module, then $i(\Lambda)$ is contained in the ball in $\LL$ (centered at the origin) with radius 
$\max\|i(\lambda_j)\|_{\LL}$. Thus, for example, if there was only one place $v|p$ of $F$ (that is, if $F=\Q$), any map
$$
LA(\rho_{\pi,\vv}) \hookrightarrow \hat{H}^0(K^p)_E[\p]
$$
automatically extends to the completion, $B(\rho_{\pi,\vv})$, giving the desired map $\eta$.

\medskip

\noindent Suppose, for simplicity, we only have two places $\{v_1,v_2\}$ of $F$ above $p$. The previous discussion gives rise to a $\GL_2(\K_{\vv_1})$-equivariant embedding,
$$
LA(\rho_{\pi,\vv_1})\hookrightarrow \Hom_{\GL_2(\K_{\vv_2})}(LA(\rho_{\pi,\vv_2}), \hat{H}^0(K^p)_E[\p]).
$$
All the homomorphisms on the right-hand side are automatically continuous. So, the target may be identified with the Banach space of bounded transformations,
$$
\LL_{\GL_2(\K_{\vv_2})}(B(\rho_{\pi,\vv_2}), \hat{H}^0(K^p)_E[\p]),
$$
which then must contain $B(\rho_{\pi,\vv_1})$. Composing with the evaluation map, yields
$$
B(\rho_{\pi,\vv_1}) \otimes B(\rho_{\pi,\vv_2}) \hookrightarrow \hat{H}^0(K^p)_E[\p],
$$
which is necessarily continuous relative to the tensor-product norm (and therefore extends to the completed tensor product): Since all invariant norms on 
$B(\rho_{\pi,\vv_i})$ are equivalent, see Corollary 5.3.4 on p. 56 in \cite{BB}, a short argument shows that {\it{any}} norm on their tensor product is dominated by (a constant multiple of) the tensor-product norm; as defined in paragraph 17 in \cite{Sc}, say.

\medskip

\noindent We may continue this, inductively, and deal with three or more places $v|p$. $\square$

\subsection{Non-vanishing at all $\p \subset \m$}

To finish the proof of Theorem 1, we have to deduce Theorem 8 from Lemma 4; again under the assumption that $\m$ is associated with a $\bar{\rho}$, with irreducible restrictions $\bar{\rho}_w$, for $w|p$. Thus, knowing that $\XX[\p]\neq 0$ for all $\p \in \CC_{\m}$, we will deduce this for all primes $\p \subset \m$ whatsoever. This proceeds exactly as the proof of Proposition 4.7 in \cite{Bre}, making use of the density of $\CC_{\m}$ (Corollary 2 in Section 4.6).

\medskip

\noindent Below we will use Nakayama's lemma, which requires the following preliminary result. 

\begin{lem}
$\Hom_{\OO}(\XX,\OO)$ is a finitely generated $\T(K^p)_{\OO,\m}$-module.
\end{lem}

\noindent {\it{Proof}}. For simplicity, let $\T=\T(K^p)_{\OO,\m}$ throughout this proof (a complete local $\OO$-algebra). By Proposition C.5 on p. 104 in \cite{Em1}, we need to show $\XX$ is cofinitely generated over $\T$ (cf. Definition C.1 in loc. cit., $\XX$ clearly satisfies the first three properties, since $\hat{H}^0(K^p)_{\OO,\m}$ is $\varpi$-adically complete, separated, and $\OO$-torsion-free). By C.1 it remains to show $(\XX/\varpi \XX)[\m]$ is finite-dimensional over $k$ (the fourth property). Note that there is a natural reduction map,
$\eta \mapsto \bar{\eta}$,
$$
\XX/\varpi\XX \longrightarrow \Hom_{\T_k[G(\Q_p)]}(\pi_{\m}/\varpi \pi_{\m}, H^0(K^p)_{k,\m}),
$$
which is injective. Since $\varpi \in \m$, after taking the $\m$-torsion, we get
$$
(\XX/\varpi\XX)[\m]\hookrightarrow \Hom_{k[G(\Q_p)]}(\otimes_{v|p}\bar{\pi}_{\m,\vv}, H^0(K^p)_{k,\m}),
$$
using that $\pi_{\m}=\otimes_{v|p}\pi_{\m,\vv}$, where $\pi_{\m,\vv}$ is a deformation of $\bar{\pi}_{\m,\vv}\leftrightarrow \bar{\rho}_{\m,\vv}$ over $\T$. We will show that the ambient space of Hom's is finite-dimensional. We use a trick from p. 78 of \cite{Em1}, from the proof of his Theorem 6.3.12. Since the representation $\bar{\pi}_{\m,\vv}$ has finite length, choose a finite-dimensional $k$-subspace $W_{\vv}$, which generates $\bar{\pi}_{\m,\vv}$  as a $\GL_2(\K_{\vv})$-representation (cf. the definition of $\VV$ in 4.4). Put $W=\otimes_{v|p}W_{\vv}$, a representation of $G(\Q_p)$ over $k$. Furthermore, since $W_{\vv}$ is smooth and finite-dimensional, we can choose a compact open subgroup $K_{\vv}$ fixing $W_{\vv}$ point-wise. Let $K_p=\prod_{v|p}K_{\vv}$. By restriction,
$$
\Hom_{k[G(\Q_p)]}(\otimes_{v|p}\bar{\pi}_{\m,\vv}, H^0(K^p)_{k,\m})\hookrightarrow \Hom_{k[K_p]}(W,H^0(K^p)_{k,\m}).
$$
Moreover, since $K_p$ acts trivially on $W$, the latter space can be thought of as 
$$
\Hom_{k[K_p]}(W,H^0(K^p)_{k,\m}) \simeq W^{\vee}\otimes_k H^0(K_pK^p)_{k,\m},
$$
which obviously has finite dimension over $k$. $\square$

\medskip

\noindent We now have everything in place to finish the proof.

\medskip

\noindent {\it{Lemma 4 $\Rightarrow$ Theorem 8}}: Let $\p \subset \m$ be a prime ideal of $\T=\T(K^p)_{\OO,\m}$. By Lemma C. 14 on p. 108 of \cite{Em1},  the natural restriction map
$$
\Hom_{\OO}(\XX,\OO)\otimes_{\T} (\T/\p) \longrightarrow \Hom_{\OO}(\XX[\p],\OO)
$$
becomes an isomorphism after tensoring $-\otimes_{\OO}E$. By the anti-equivalence of Proposition C.5 in loc. cit., it therefore suffices to show that
$$
\text{$M/\p M\neq 0$, $\y$ $M=\Hom_{\OO}(\XX,\OO)$.} 
$$
Once we show $\T$ acts faithfully on $M$, we are done by Nakayama's lemma (which applies since $M$ is finitely generated, as shown above): Suppose $M=\p M$. Then there is a $t \equiv 1$ (mod $\p$) in $\T$ such that $tM=0$. Clearly a contradiction. 

\medskip

\noindent To show faithfulness, first note that $M/\p M$ is a vector space over $\T/\p$, a finite field extension of $E$. Thus, $\T/\p$ acts faithfully on $M/ \p M$, whenever the latter is nonzero. If $t \in \T$ acts trivially on $M$, it acts trivially on every $M/ \p M$, and therefore $t$ belongs to every $\p$ for which 
$M/ \p M$ is nonzero. If this holds for a Zariski dense set $\mathcal{S}$ of primes, take $\mathcal{S}=\CC_{\m}$ for instance, we infer that
$$
t \in \cap_{\p \in \mathcal{S}} \p=0.
$$ 
That is, $t=0$; which is to say $\T$ acts faithfully on $M$. $\square$



\noindent {\it{Przemyslaw Chojecki}}

\noindent {\sc{Institut Mathematique de Jussieu, Paris, France.}}

\noindent {\it{E-mail address}}: {\texttt{chojecki@math.jussieu.fr}}

\bigskip

\noindent {\it{Claus Sorensen}}

\noindent {\sc{Department of Mathematics, Princeton University, USA.}}

\noindent {\it{E-mail address}}: {\texttt{claus@princeton.edu}}


\begin{thebibliography}{}

\bibitem[BLGGT1]{B1} T. Barnet-Lamb, T. Gee, D. Geraghty and R. Taylor,
{\it{Local-global compatibility for l=p. I}}, Ann. de Math. de Toulouse 21 (2012), 57-92.

\bibitem[BLGGT2]{B2} T. Barnet-Lamb, T. Gee, D. Geraghty and R. Taylor, {\it{Local-global compatibility for l=p. II }}, to appear Ann. Sci. de l'ENS. 

\bibitem[BC]{BC} J. Bellaiche and G. Chenevier, {\it{Families of Galois representations and Selmer groups}}. Asterisque 324 (2009), Soc. Math. France. 

\bibitem[Ber]{Ber} L. Berger, {\it{La correspondance de Langlands locale p-adique pour $\GL_2(\Q _p)$}}. Asterisque No. 339 (2011), 157-180.

\bibitem[BB]{BB} L. Berger and C. Breuil, {\it{Sur quelques representations potentiellement cristallines de $\GL_2(\Q _p)$}}. Asterisque 330 (2010), 155-211.

\bibitem[BE]{BE} C. Breuil and M. Emerton, {\it{Representations p-adiques ordinaires de $\GL_2(\Q_p)$ et compatibilite local-global}}. Asterisque 331 (2010) 255-315.

\bibitem[BS]{BS} C. Breuil and P. Schneider, {\it{First steps towards p-adic Langlands functoriality}}. J. Reine Angew. Math. 610, 2007, 149-180.


\bibitem[Bre]{Bre} C. Breuil, {\it{Correspondance de Langlands p-adique, compatibilite local-global et applications}}. Bourbaki seminar, January 2011, to appear in Asterisque. 

\bibitem[Buz]{Buz} K. Buzzard, {\it{Eigenvarieties}}. In "L-functions and Galois representations", 59-120, London Math. Soc. Lecture Note Ser., 320, Cambridge Univ. Press, Cambridge (2007).

\bibitem[Car]{Car} A. Caraiani, {\it{Monodromy and local-global compatibility for l=p}}. Preprint (2012).

\bibitem[Che1]{Ch} G. Chenevier, {\it{Sur la densite des representations cristallines du groupe de Galois absolu de $\Q _p$}}. Math. Annalen 335, 1469-1525 (2013).

\bibitem[Che2]{Che} G. Chenevier, {\it{Une application des varietes de Hecke des groupes unitaires}}. Preprint (2009).

\bibitem[Cho]{Cho} P. Chojecki, {\it{Density of crystalline points on unitary Shimura varieties}}. Int. J. Number Theory, Vol. 9, No. 3 (2013) 1-15.

\bibitem[CM]{CM} R. Coleman and B. Mazur, {\it{The eigencurve}}. In "Galois representations in arithmetic algebraic geometry (Durham, 1996)", 1-113,
London Math. Soc. Lecture Note Ser., 254, Cambridge Univ. Press, Cambridge (1998).

\bibitem[Co1]{Col} P. Colmez, {\it{Representations de $\GL_2(\Q _p)$ et $(\phi,\Gamma)$-modules}}. Asterisque 330 (2010), 281-509.

\bibitem[Co2]{Co} P. Colmez, {\it{La serie principale unitarire de $\GL_2(\Q_p)$: Vecteurs localement analytiques}}. Preprint (2012).

\bibitem[CY]{CY} T. Barnet-Lamb, D. Geraghty, M. Harris, and R. Taylor, {\it{A family of Calabi-Yau varieties and potential automorphy, II}}. P.R.I.M.S. 47 (2011), 29-98.

\bibitem[Em1]{Em0} M. Emerton, {\it{A local-global compatibility conjecture in the p-adic Langlands programme for $\GL _{2 / \Q}$ }}. Pure and Applied Math. Quarterly 2 no. 2, 279-393 (2006).

\bibitem[Em2]{Em1} M. Emerton, {\it{Local-global compatibility in the p-adic Langlands programme for $\GL_{2 / \Q }$}}. Preprint (2011).

\bibitem[Em3]{Em2} M. Emerton, {\it{Locally analytic vectors in representations of locally p-adic analytic groups}}. To appear in Memoirs of the AMS.

\bibitem[Em4]{Em} M. Emerton, {\it{On the interpolation of systems of eigenvalues attached to automorphic Hecke eigenforms}}. Invent. Math. 164, no. 1, 1-84 (2006).

\bibitem[EH]{EH} M. Emerton and D. Helm, {\it{The local Langlands correspondence for $\GL(n)$ in families}}. Preprint (2012).

\bibitem[HT]{HT} M. Harris and R. Taylor, {\it{The geometry and cohomology of some simple Shimura varieties}}. Annals of Math. Studies 151, PUP (2001).

\bibitem[He]{He} E. Hellmann, {\it{Families of trianguline representations and finite slope spaces}}. Preprint (2012).

\bibitem[Ki]{Ki} M. Kisin, {\it{Deformations of $G_{\Q _p}$ and $\GL_2(\Q_p)$ representations}}. Asterisque 330, (2010), 529-542.

\bibitem[Loe]{Loe} D. Loeffler, {\it{Overconvergent algebraic automorphic forms}}. Proc. London Math. Soc. 102 (2011), no. 2, 193-228.

\bibitem[LXZ]{LXZ} R. Liu, B. Xie, and Y. Zhang, {\it{Locally analytic vectors of unitary principal series of $\GL_2(\Q_p)$}}. Preprint (2011).

\bibitem[Min]{Min} A. Minguez, {\it{Unramified representations of unitary groups}}. In  “Stabilisation de la formule des traces, varietes de Shimura, et applications arithmetiques”, Int. Press of Boston.

\bibitem[Pas]{Pas} V. Paskunas, {\it{On the image of Colmez's Montreal functor}}. To appear in Publications mathematiques de l'IHES.

\bibitem[Rog]{Rog} J. Rogawski, {\it{Automorphic representations of unitary groups in three variables}}. Annals of Mathematics Studies, 123, Princeton University Press, Princeton, NJ (1990).

\bibitem[Sc]{Sc} P. Schneider, {\it{Nonarchimedean functional analysis}}. Springer Monographs in Mathematics, Springer-Verlag, Berlin (2002).

\bibitem[So1]{S} C. Sorensen, {\it{A proof of the Breuil-Schneider conjecture in the indecomposable case}}. Annals of Mathematics 177 (2013), 1-16.  

\bibitem[So2]{Sor} C. Sorensen, {\it{Eigenvarieties and invariant norms: Towards p-adic Langlands for U(n)}}. Preprint (2012).




\end{thebibliography}
\end{document}